\documentclass{article}

\usepackage[utf8]{inputenc}
\usepackage[top=2.5cm, bottom=2.5cm, left=2.5cm, right=2.5cm]{geometry}

\usepackage{amssymb,amsmath}
\usepackage{mdwlist} 
\usepackage{graphicx} 
\usepackage{bm} 
\usepackage{subfig} 	

\usepackage{placeins}

\usepackage[framemethod=TikZ]{mdframed}

\usepackage{graphicx}

\usepackage{textcomp} 

\usepackage{xcolor}

\title{Problem-Driven Scenario Clustering in Stochastic Optimization}
\author{Julien Keutchayan, Janosch Ortmann, Walter Rei}
\date{\today}

\begin{document}

\maketitle

\newcommand{\janosch}[1]{(\textbf{\textcolor{blue}{#1}})}
\newcommand{\walter}[1]{(\textbf{\textcolor{red}{#1}})}
\newcommand{\julien}[1]{(\textbf{\textcolor{green}{#1}})}

\newcommand{\E}{\mathbb{E}}
\newcommand{\R}{\mathbb{R}}
\newcommand{\N}{\mathbb{N}}
\newcommand{\whx}{\widehat{x}}
\newcommand{\wtx}{\widetilde{x}}
\newcommand{\argmin}{\mathrm{argmin}}
\newcommand{\rt}{\textsuperscript \textregistered}

\newcommand{\f}{\mathcal{F}}
\newcommand{\y}{\mathcal{Y}}
\newcommand{\ee}{\mathbb{E}}
\newcommand{\pp}{\mathbb{P}}
\newcommand{\rr}{\mathbb{R}}

\newcommand{\abs}[1]{\left|#1\right|}
\newcommand{\ab}[1]{\left[#1\right]}
\newcommand{\rb}[1]{\left(#1\right)}
\newcommand{\set}[1]{\left\{#1\right\}}

\newcommand{\setz}[2]{\set{#1,\dots,#2}}

\newcommand{\twovec}[2]{\begin{pmatrix} #1 \\ #2 \end{pmatrix}}

\begin{abstract}
In stochastic optimisation, the large number of scenarios required to faithfully represent the underlying uncertainty is often a barrier to finding efficient numerical solutions. This motivates the \emph{scenario reduction problem}: by find a smaller subset of scenarios, reduce the numerical complexity while keeping the error at an acceptable level.

In this paper we propose a novel and computationally efficient methodology to tackle the scenario reduction problem when the error to be minimised is the \emph{implementation error}, i.e. the error incurred by implementing the solution of the reduced problem in the original problem. Specifically, we develop a problem-driven scenario clustering method that produces a partition of the scenario set. Each cluster contains  a representative scenario that best reflects the conditional objective values in each cluster of the partition to be identified.

We demonstrate the efficiency of our method by applying it to two challenging stochastic combinatorial optimization problems: the two-stage stochastic network design problem and the two-stage facility location problem. When compared to alternative clustering methods and Monte Carlo sampling, our method is shown to clearly outperform all other methods.
\end{abstract}

\section{Introduction}
\label{Section_Introduction}

When uncertainty is present in decision-making contexts, scenario generation refers to the methods used to sample the set of possible outcomes of the random parameters, see \cite{kaut2012scenario, keutchayan2017quality}.
Such methods are staples in the design of decision support techniques that explicitly consider the uncertainty affecting the parameters defining the informational context in which problems are to be solved.
In particular, scenario generation methods are readily applied in stochastic optimization to obtain approximations of the probability distributions, i.e., samples of possible outcomes, that are used to implicitly formulate how the uncertain (or stochastic) parameters randomly vary, see \cite{birge2011introduction}.
The generation of such samples can then be used statically to solve particularly hard stochastic optimization models, such as the stochastic network design problems (see, \cite{Rahmaniani:SJOO:2018}); or probabilistically to produce probabilistic bounds, such as in the sample average approximation method (e.g., \cite{Kleywegt:SJOO:2002,Luedtke:SJOO:2008}) and the stochastic decomposition strategy (e.g., \cite{Higle:MOOR:1991}).
In all cases, an important challenge that is often faced is that the size of the generated scenario set can be very large, thus hindering its use for the computational purposes associated with stochastic optimization.

Specifically, when a scenario set is used to formulate a stochastic optimization model, the resulting scenario-based formulation will typically be large-scale with its size being directly related to the cardinality of the set, see e.g. \cite{birge2011introduction}. Such formulations quickly become intractable to solve directly, especially in the case of stochastic integer programs, as shown for instance in \cite{Dyer:MP:2006}. Therefore, the generation of scenarios is a delicate step in any algorithmic strategy aimed to efficiently solve stochastic combinatorial optimization models. This step is subject to a trade-off: on the one hand, a sufficiently large set of meaningful scenarios should be generated to properly capture how uncertain parameters may randomly vary; but on the other hand, the size of the set should not be too large that it prevents the resulting formulation to be solved in an appropriate manner.

The combinatorial challenges of large-scale optimization models motivate the study of \emph{scenario reduction} methods. The goal of these methods can be stated as follows: given a particular scenario set, identify a subset of a given fixed size that minimizes some approximation error induced by replacing the original set with the identified subset. Considering that the number of scenarios required to model appropriately a real-word (generally multi-dimensional) probability distribution is typically too large to be directly used in the context of solving stochastic combinatorial optimization problems, tackling the scenario reduction problem is unavoidable to design computational efficient solution approaches in this context. 


When surveying the scientific literature, one finds that there has been a steady stream of research dedicated to the development of methods to efficiently generate scenario sets that can be applied in the context of performing stochastic optimization. Such methods can be classified in two groups: the \emph{distribution-driven} ones and the \emph{problem-driven} ones. Distribution-driven methods perform scenario generation based on a criterion linked to the probability distribution that underlies the problem, and not directly the problem itself. Such approaches seek to find the subset that provides the approximation that is probabilistically closest to the original scenario set, see, e.g., \cite{Dupacova:MP:2003, Romisch:Proceedings:2009, Henrion:COA:2009}. In this case, a metric is used to evaluate the distance between probability distributions, which defines the error function that is minimized by the scenario reduction procedure. This is the case for instance of the methods that minimize the Wasserstein distances (e.g., \cite{pflug2015dynamic}), or those that aim at matching the moments of the distribution (e.g., \cite{hoyland2003heuristic}). Other scenario generation methods such as Monte Carlo (e.g., \cite{Kleywegt:SJOO:2002}), quasi-Monte Carlo methods (e.g., \cite{heitsch2016quasi}), and sparse grid quadrature (e.g., \cite{chen2015scenario}) also fall in this category since they aim at discretizing the probability distribution using a finite number of samples. Stated differently, a distribution-driven method will tend to operate in the same way when applied to two different stochastic problems with the same underlying probability distribution. 

Conversely, a problem-driven method aims at generating scenarios by considering the problem as a whole, which in addition to the probability distribution also includes any relevant property of the problem that can be drawn from its objective function and constraints. In that case, two problems with the same underlying distribution but different objective functions will likely not be formulated using the same scenario set. The objective function is indeed important to consider since, for instance, a risk-averse problem featuring a log-utility function might require a different scenario set than its risk-neutral version as it might be sensitive to different aspects of the distribution. Although the vast majority of methods used for scenario generation fall in the distribution-driven category, as it has been the only predominant one historically, approaches falling in the problem-driven category have been recently studied in \cite{Henrion:MP:2018, fairbrother2019problem, Keutchayan:MOOR:2020} as a mean to provide more efficient procedures to generate scenarios.

The problem-driven approaches are logically viewed as being more powerful, since they are made to consider the problem as a whole and not solely its probability distribution. However, research in this direction has been very sparse and essentially limited to theoretical results. Practitioners currently lack practical algorithms that could be used systematically to generate scenarios for their problem. The situation is different for the distribution-driven approaches, which provide plenty of algorithms available off the shelf, such as Monte Carlo sampling or the $k$-means clustering. The goal of our paper is to provide one such systematic problem-driven algorithm.

We propose a novel and computationally efficient methodology to tackle the scenario reduction problem when the objective pursued is to minimize the implementation error. This error is defined as the difference between the optimal value of the original problem and the value of the solution provided by solving the problem on the reduced scenario set. In practice, it corresponds to the error that would be made by implementing the scenario tree decisions in the original problem. It is therefore very relevant for practitioners. Specifically, we develop a problem-driven scenario clustering method that produces a partition of the scenario set that enables representative scenarios to be identified, which best reflect the conditional objective values over the cluster they represent.

The proposed method is applied to solve the scenario reduction problem for two challenging stochastic combinatorial optimization problems: the two-stage stochastic network design problem and the two-stage facility location problem.
When compared to four alternative methods:  $k$-means, $k$-medoids, $k$-medians and Monte Carlo sampling, our method is shown to clearly outperform all other methods, i.e., it greatly reduces the number of scenarios while efficiently limiting the implementation errors induced.
Moreover, in the case of the two-stage facility location problem, the method is capable of finding near-optimal solutions to the full size scenario-based model by reducing the original set to a single scenario, which the other tested methods are incapable of.

The rest of the paper is structured as follows. Section \ref{sec:problem} presents the general stochastic optimization problem that we are considering in this paper. Section~\ref{scenario_reduction_methods} describes some of the most popular methods that can be used to reduce the number of scenarios. Section \ref{sec:algo} is dedicated to the statement of our algorithm along with a concrete illustration. Section \ref{sec:numerical} presents our numerical results, and Section \ref{sec:conclusion} concludes.

\section{Problem Statement}
\label{sec:problem}

This section is devoted to a precise statement of the class of stochastic optimization problems for which we develop a new problem-driven scenario reduction method. 

\subsection{Original Problem}
\label{sec_original_problem}
We consider the following stochastic optimization problem:

\begin{equation}
\label{eq_original_problem}
\min_{x\in X\subseteq \R^n} \left\{f(x) := \frac{1}{N} \sum_{i=1}^N F(x, \xi_i)\right\},
\end{equation}
where $\xi_1,\dots,\xi_N$ are equiprobable scenarios taking values in $\R^d$ and $F(x, \xi_i)$ represents the cost associated to the decisions $x\in X$ in scenario $\xi_i$. The optimal solution of this problem is denoted by $x^*$ and its optimal value by $v^*$. Note that there is no restriction in the number of scenarios ($N$), decisions ($n$), and random parameters ($d$); these can be any positive integers. Note also that the problem is formulated using equiprobable scenarios for simplicity; the method introduced in this paper can be easily generalized to scenarios with different probabilities. 

The cost function $F$ may be given explicitly (if the problem is one-stage), or may be itself the result of a second-stage optimization problem. In this latter case, which is the framework of two-stage stochastic programming, it typically takes the following form:
\begin{equation}
\label{eq_2nd_stage_problem}
F(x,\xi_i) = \min_{y\in Y(x,\xi_i)} g(x, y, \xi_i),
\end{equation}
for some second-stage cost function $g(x, y, \xi_i)$ and second-stage decisions $y\in Y(x,\xi_i)\subseteq \R^{m}$ (see, for instance, \cite{birge2011introduction}). 

In the following, we require the problem to satisfy the condition of \emph{relative complete recourse}. That is, it must hold that $F(x,\xi_i) <\infty$ for all $x\in X$ and $i\in\{1,\dots,N\}$. The reason this property is required will become clear in the next section when we introduce the approximate problem. We do not impose any other condition on $F$ besides the basic ones ensuring that problem \eqref{eq_original_problem} is well-defined, i.e., there exists an optimal solution $x^*$. Such conditions can be found in \cite{rockafellar2009variational}. Note that, in practice, some conditions like the linearity (or convexity) of $F(x,\xi_i)$ as a function of $x$ greatly improve one's ability to solve the problem; however they are not technically required in our theoretical developments.

The scenarios $\xi_1,\dots,\xi_N$ may be given from historical data (e.g., from $N$ daily observations of demands for $d$ items in a retail store), or they may be sampled from a multivariate probability distribution as in the sample average approximation (SAA) of \cite{shapiro2003monte}. 
In some cases, if $N$ is large (how large depends on the problem) and if the problem is computationally challenging to solve, it may happen that finding a good solution for problem \eqref{eq_original_problem} with $N$ scenarios is out of reach for even the best solvers available today. In that situation, one way to ease the computational burden is to reduce the problem size by building a smaller problem with fewer scenarios. However, building such a problem is not an easy task as it should be done with the goal to find a solution with a value as close as possible to that of the original problem, $v^*$.

\subsection{Approximate Problem}

Suppose that problem \eqref{eq_original_problem} cannot be solved as it is, so we need to build from it an approximate problem composed of $K$ scenarios with $K \ll N$. There are two broad ways to generate those scenarios: they may be picked directly in the original set $\{\xi_1,\dots,\xi_N\}$, or they may be completely new scenarios that do not exist in the original set. In any case, these scenarios are generated for the purpose of defining an approximate problem close to the original one in a sense that we will make precise below. 

Consider for now that this set has been computed, and let us denoted it by $\{\widetilde{\xi}_1, \dots, \widetilde{\xi}_K\}$ with the corresponding probabilities $\{p_1,\dots,p_K\}$. The approximate problem takes the form:
\begin{equation}
\label{eq_approximate_problem}
\min_{x\in X\subseteq \R^n} \left\{\widetilde{f}(x) := \sum_{k=1}^K p_k F(x, \widetilde{\xi}_k) \right\},
\end{equation}
and its optimal solution is denoted by $\wtx^*$.

Under the condition of relative complete recourse stated in Section~\ref{sec_original_problem}, any optimal solution of \eqref{eq_approximate_problem} is a feasible solution of \eqref{eq_original_problem}. This ensures that $\wtx^*$ can be evaluated in the original problem. Its value, $f(\wtx^*)$, therefore provides an upper bound on $v^*$. 

The goal of the scenario reduction approach is to find the scenarios $\widetilde{\xi}_1, \dots, \widetilde{\xi}_K$ and their probabilities $p_1, \dots, p_K$ that provide the smallest upper bound, i.e., such that the gap
\begin{equation}
\label{implementation_error}
f(\wtx^*) - v^* = \frac{1}{N} \sum_{i=1}^N F(\wtx^*, \xi_i) - \frac{1}{N} \sum_{i=1}^N F(x^*, \xi_i) \geq 0,
\end{equation}
is as close as possible to zero. 

In the following, we refer to \eqref{implementation_error} as the \emph{implementation error}, as it corresponds to the error that would be made by implementing in the original problem \eqref{eq_original_problem} the optimal decisions given by the approximate problem \eqref{eq_approximate_problem}. It is always non-negative (by optimality of $x^*$) and it equals zero if and only if $\wtx^*$ is also optimal for \eqref{eq_original_problem}. 

\section{Scenario Reduction Methods}
\label{scenario_reduction_methods}

Clustering methods provide a generic way to compute the reduced set of scenarios $\{\widetilde{\xi}_1, \dots, \widetilde{\xi}_K\}$. Generally speaking, a clustering method is an algorithmic procedure that partitions the original scenario set $\{\xi_1,\dots,\xi_N\}$ into $K$ clusters $C_1,\dots,C_K$, of possibly different sizes, such that scenarios close to each other in the sense of some metric $d$ lie in the same cluster. In the context of scenario reduction, we are specifically interested in the clustering methods that output a centroid scenario $\widetilde{\xi}_k$ for each cluster $C_k$; this centroid becomes the representative scenario of the cluster and is then given the probability $p_k = \vert C_k\vert/ N$. Clustering methods will differ by the metric used to measure the distance between each pair of scenarios and by whether they pick the representatives as a member of the original set or not, i.e., whether $\widetilde{\xi}_k \in C_k$ or $\widetilde{\xi}_k \not\in C_k$. For a general presentation on clustering methods, see, e.g., \cite{han2011data, jain1988algorithms}.

Perhaps the most popular clustering method is the \emph{$k$-means} algorithm, which aims at minimizing the sum of square $L_2$-distances within each cluster (see, e.g., \cite{ilprints778}). The centroid is given by the average scenario over the cluster; as such it may not be a member of the original set. The $k$-means algorithm is a method of vector quantization particularly suitable for scenario reduction. It is connected to the problem of finding the distribution sitting on $K$ scenarios that minimizes the 2nd-order Wasserstein distance from the original scenarios set (see \cite{ho2017multilevel}). This minimization problem is computationally difficult (NP-hard). The $k$-means algorithm provides a fast heuristic approach that allows finding locally optimal solutions quickly. For more information on the Wasserstein distance and its relevance in the context of stochastic optimization, we refer to \cite{pflug2001scenario, pflug2015dynamic}. 

The \emph{$k$-medians} algorithm (\cite{jain1988algorithms}) follows a similar idea but minimizes the sum of $L_1$-distances. It has the advantage of being less sensitive to outliers than $k$-means, because of the robustness of the median as compared to the mean. The $k$-means and $k$-medians algorithms are suitable for numerical data but not for scenarios in the form of categories. The \emph{$k$-modes} algorithm (\cite{huang1998extensions}) considers the discrete distance, which equals zero if the two scenarios are equal and one otherwise, which makes it suitable to work with scenarios made of categorical variables. 

The three aforementioned methods have in common that the representative scenarios may not lie in the original set. The \emph{$k$-medoids} algorithm (also sometimes called ``partition around medoids'') is similar to $k$-means, but unlike the latter, it enforces the centroids to be part of the original set of scenarios. Moreover, it can be used with any scenario distance $d$ (see, e.g., \cite{schubert2019faster}). In the context of stochastic optimization, the $k$-medoids algorithm provides a local search heuristic that aims at solving the optimal scenario reduction problem formulated in \cite{rujeerapaiboon2018scenario, Henrion:MP:2018}. A practical advantage of $k$-medoids over $k$-means is that the former guarantees that the scenarios used in the approximate problem are actual scenarios that can occur in reality since they are picked in the original set. However, the latter has the advantage to be much faster.

Recently, another approach that has been developed to group scenarios that exhibit common characteristics is the decision-based clustering method introduced in \cite{Hewitt:2021aa}.
This approach seeks to cluster scenarios (from the original set) that induce similar solutions with respect to the stochastic optimization model.
For a given scenario, its induced solution represents the decisions that would be made if the particular scenario could be accurately predicted.
Based on this clustering approach, an alternative medoid upper bound was proposed to obtain an approximate problem.
In this case, the medoid of a cluster is identified as the scenario whose induced solution represents the closest proxy for the solutions associated to all other scenarios within the group.

Finally, scenario reduction can also be performed using Monte Carlo sub-sampling. Unlike clustering methods, which must solve an optimization problem and therefore are computationally intensive, Monte Carlo sub-sampling simply consists in picking randomly a subset of $K$ scenarios in the original set and assigning them equal probability $1/K$. The relevance of Monte Carlo methods in the context of stochastic optimization has been well studied in \cite{shapiro2000rate, shapiro2003monte}.

\section{A New Scenario Clustering Approach}
\label{sec:algo}

In this section, we introduce the \emph{cost-space scenario clustering} (CSSC) algorithm and illustrate its benefit using a concrete example. The development of a clustering method driven by the cost function of the problem is first motivated in subsection~\ref{motivation_cost_driven}. The mathematical basis is then provided in subsection~\ref{mathematical_developments} and the algorithm is detailed in subsection~\ref{Section_algo_description}. In subsection~\ref{Section_method_illustration}, the algorithm is applied on a specific two-stage stochastic problem to both illustrate its use and contrast its results with those obtained by applying clustering methods purely driven by the distribution.

\subsection{Motivation of Cost-Driven Clustering}
\label{motivation_cost_driven}

The clustering approach developed in this paper does not consider a distance between scenarios in the \emph{space of random outcomes} ($\R^d$) but instead in the \emph{space of cost values} ($\R$). 
Specifically, this means that the proximity between two scenarios $\xi_i$ and $\xi_j$ is never measured directly by some metric $d(\xi_i, \xi_j)$ defined over $\R^d\times \R^d$ but through their value in the cost function: $F(\cdot,\xi_i)$ and $F(\cdot,\xi_j)$. This choice is motivated by the fact that, ultimately, a scenario reduction approach should aim at finding representative scenarios that provide a good approximation of the original problem, regardless of whether the clusters turn out to be good or bad in the sense of any such metric $d$. 

To better understand this, consider the case where two scenarios $\xi_i$ and $\xi_j$ are very different, i.e., $d(\xi_i, \xi_j)$ is large, but where they lead to the same value of the cost function, i.e., $F(\cdot,\xi_i) = F(\cdot, \xi_j)$. It is then easy to see that these two scenarios can be grouped together into the same cluster with either representative $\xi_i$ or $\xi_j$ and weight $\tfrac2N$ without changing the solution of the original problem, hence resulting in a reduction of scenarios without loss of information as far as the stochastic problem is concerned. 

The previous case of pairwise equality of cost values is not the only one that allow an error-free reduction. Consider now three scenarios $\xi_i$, $\xi_j$, $\xi_k$ such that the cost function of $\xi_k$ equals the mean of the other two, i.e., $F(\cdot,\xi_k) = \tfrac12 (F(\cdot, \xi_i) + F(\cdot, \xi_j))$. It is then possible to cluster together $\set{\xi_i,\xi_j,\xi_k}$ with weight $\tfrac3N$ and representative $\xi_k$ without changing the solution of the original problem, hence effectively removing two-thirds of the scenarios without introducing any error. It is easy to see that this situation can be generalized to any number of scenarios, as long as there is one scenario (the representative one) with a cost function that can be expressed as the mean of the other ones. 

It is clear that if these cases arise in practice, it will be with an approximate equality rather than an exact one, which means that some error is likely to be introduced by the clustering. However, if one has the goal to keep this error as low as possible, this motivates the development of a reduction method that clusters scenarios in the space of cost values to spot such similarities that are not apparent to distribution-based methods. The next sections will provide further explanations and illustrate this on a concrete stochastic problem.

\subsection{Mathematical Developments}
\label{mathematical_developments}

In this section, we derive the error bound to be minimized by the CSSC algorithm from the implementation error defined in \eqref{implementation_error}. To this end, let us first decompose the implementation error as follows:
\begin{align}
\label{implementation_error_decomposition_1}
\vert f(\wtx^*) - v^* \vert = \vert f(\wtx^*) - f(x^*) \vert & = \vert f(\wtx^*) - \widetilde{f}(\wtx^*) + \widetilde{f}(\wtx^*) - f(x^*) \vert \\
\label{implementation_error_decomposition_2}
& \leq \vert f(\wtx^*) - \widetilde{f}(\wtx^*) \vert + \vert \widetilde{f}(\wtx^*) - f(x^*) \vert.
\end{align}
The second term can be further bounded by:
\begin{align}
\label{implementation_error_decomposition_3}
\vert \widetilde{f}(\wtx^*) - f(x^*) \vert & = \max\{\widetilde{f}(\wtx^*) - f(x^*), \, f(x^*) - \widetilde{f}(\wtx^*)\} \\
\label{implementation_error_decomposition_4}
& \leq \max\{\widetilde{f}(x^*) - f(x^*), \, f(\wtx^*) - \widetilde{f}(\wtx^*)\} \\
\label{implementation_error_decomposition_5}
& \leq \max\{\vert \widetilde{f}(x^*) - f(x^*) \vert, \, \vert f(\wtx^*) - \widetilde{f}(\wtx^*) \vert\} \\
\label{implementation_error_decomposition_6}
& = \max_{x\in\{x^*, \wtx^*\}} \vert \widetilde{f}(x) - f(x) \vert.
\end{align}
where for \eqref{implementation_error_decomposition_4} we use the fact that $\widetilde{f}(\wtx^*) \leq \widetilde{f}(x^*)$ since $\wtx^*$ is a minimum of $\widetilde{f}$, and $f(x^*)\leq f(\wtx^*)$ since $x^*$ is a minimum of $f$. 

Finally, by combining \eqref{implementation_error_decomposition_2} and \eqref{implementation_error_decomposition_6} we can bound the implementation error as follows. Let $\widetilde{X}\subseteq X$ be any feasible set such that $\{x^*, \wtx^*\}\subset\widetilde{X} $, that is $\widetilde X$ includes the two optimal solutions of the original and approximate problems. 
\begin{align}
	\vert f(\wtx^*) - v^* \vert & \leq 2 \max_{x\in\{x^*, \wtx^*\}} \vert \widetilde{f}(x) - f(x) \vert \\
	& \leq 2 \, \max_{x\in \widetilde{X}} \, \vert \widetilde{f}(x) - f(x) \vert \\
	& = 2\max_{x\in \widetilde{X}} \left\vert \frac{1}{N} \sum_{i=1}^N F(x, \xi_i) -\sum_{k=1}^K p_k F(x, \widetilde{\xi}_k) \right\vert \\
	& = 2\max_{x\in \widetilde{X}} \left\vert \frac{1}{N} \sum_{k=1}^K \sum_{i\in C_k} F(x,\xi_i) - \sum_{k=1}^K \frac{\vert C_k \vert}{N} F(x, \widetilde{\xi}_k) \right\vert \\
	& = 2\max_{x\in \widetilde{X}} \left\vert \sum_{k=1}^K \frac{\vert C_k \vert}{N} \left( \frac{1}{\vert C_k \vert} \sum_{i\in C_k} F(x,\xi_i) - F(x,\widetilde{\xi}_k) \right) \right\vert \\
	& \leq 2\sum_{k=1}^K p_k \sup_{x\in \widetilde{X}} \left\vert \frac{1}{\vert C_k \vert} \sum_{i\in C_k} F(x,\xi_i) - F(x,\widetilde{\xi}_k) \right\vert \\
	\label{clustering_discrepancy}
	& =: 2\sum_{k=1}^K p_k D(C_k).
\end{align}

The quantity $D(C_k)$ can be seen as the \emph{discrepancy} of the cluster $C_k$. It measures how much the cost function $F(x,\widetilde{\xi}_k)$ of its representative scenario $\widetilde{\xi}_k$ matches the average cost values of the whole cluster $C_k$ over the feasible set $\widetilde{X}$. It should be noted that this discrepancy generalizes mathematically to any number of scenarios the cost-driven clustering idea that was expressed in subsection~\ref{motivation_cost_driven} using the example of three scenarios.

A true measure of discrepancy (in the sense of theoretically valid) would have to consider a set $\widetilde{X}$ sufficiently large to ensure the validity of the bound, while sufficiently small to guarantee the tightness of the bound. In practice, however, this trade-off will be difficult to resolve. The set $\widetilde{X}$ would have to include the optimal solutions of the original and approximate problems, which are both unknown at the time of the scenario reduction. A way to bypass this issue would be to consider a larger set that includes any ``sensible'' solutions, or simply setting $\widetilde{X}=X$ if one cannot define such a set. However, this would lead to a large set of solutions (possibly infinitely many), which would render the computation of $D(C_k)$ nearly impossible in practice. Note that computing $F(x,\xi_i)$ for a single $x$ requires the solution of the 2nd-stage problem, hence it is not free of computational cost. 

In order to obtain a usable bound, we propose to consider a set $\widetilde{X}$ composed of all the optimal solutions of the one-scenario subproblems. 
These solutions are easier to compute, as they require solving the deterministic version of the problem composed of a single scenario, and together they provide a proxy of sensible solutions.
Furthermore, it was illustrated in \cite{Hewitt:2021aa} that these solutions could be used to effectively search for commonalities between the scenarios.
We leverage them here to define the basis on which the discrepancy measure is evaluated.
Thus, we approximate $D(C_k)$ by: 
\begin{equation}
\label{disrepancy_approx}
	D(C_k) \simeq \left\vert \frac{1}{\vert C_k \vert} \sum_{i\in C_k} F(x_k^*,\xi_i) - F(x_k^*,\widetilde{\xi}_k) \right\vert,
\end{equation}
where $x_k^*$ is the optimal solution of the deterministic version of problem \eqref{eq_original_problem} with scenario $\widetilde{\xi}_k$, i.e., $x_k^* \in \underset{x\in X}{\argmin}\;F(x,\widetilde{\xi}_k)$.

The next subsection describes the algorithm in more details and analyze its computational cost.

\subsection{Algorithm Description}\label{Section_algo_description}

The CSSC algorithm aims to minimize the clustering bound defined in \eqref{clustering_discrepancy} using the discrepancy approximation \eqref{disrepancy_approx}. Overall, it consists in the following two steps:

\begin{itemize}
	\item \textbf{Step 1.} Compute the opportunity-cost matrix $\bm{V} = (V_{i,j})$ where 
	\begin{equation}
	\label{eq_cost_value_matrix}
	V_{i,j} = F(x_i^*,\xi_j), \qquad \forall i,j\in\setz1N,
	\end{equation}
	where each $x_i^*$ is the optimal solution of the one-scenario subproblem:
	\begin{equation}
	\label{eq_one_scenario_problem}
	x_i^* \in \underset{x\in X}{\argmin}\;F(x,\xi_i), \qquad \forall i\in\setz1N.
	\end{equation}
	\item \textbf{Step 2.} Find a partition of the set $\setz1N$ into $K$ clusters $C_1,\dots,C_K$ and their representatives $r_1\in C_1,\dots,r_K\in C_K$ such that the following \emph{clustering discrepancy} is minimized:
	\begin{equation}
	\label{clustering_error}
	\sum_{k=1}^K p_k \left\vert V_{r_k,r_k} - \frac{1}{\vert C_k\vert} \sum_{j\in C_k} V_{r_k,j}\right\vert,
	\end{equation}
	where $p_k = \frac{\vert C_k\vert}{N}$.
\end{itemize}

The first step of the algorithm consists in solving $N$ times the one-scenario subproblem \eqref{eq_one_scenario_problem}, which becomes a deterministic problem, and then evaluating the resulting solutions $x_i^*$, $i \in\setz1N$, with respect to each scenario $\xi_j$, i.e., thus computing the values $V_{i,j} = F(x_i^*,\xi_j)$, $\forall i,j \in\setz1N$. 
In the case where the original problem is two-stage, these evaluations require solving the second-stage problem $N^2$ times (or more precisely $N^2 - N$ times, since the values $V_{i,i}$, for $i \in\setz1N$, are already given by \eqref{eq_one_scenario_problem}). 
The computational effort to perform these nearly $N^2$ iterations may seem prohibitive if the problem has integer decisions.
However, as we will see in the two case studies in subsections \ref{network_design} and \ref{facility_location}, this can be easily worked around by considering the linear relaxation of the problem. 
Moreover, it should be noted that the scenario clustering method that is developed seeks to group scenarios that exhibit commonalities defined in terms of the problem in which the scenarios are used.
Of course, such commonalities can be searched for by directly considering the solutions and their associated costs as defined in the original problem, i.e., the evaluations as stated by (\ref{eq_cost_value_matrix}) and (\ref{eq_one_scenario_problem}).
However, this search can also be performed by considering any simplification of the original problem (e.g., either a relaxation or a restriction).
If a problem simplification is used to perform the cluster search then the resulting clustering would necessarily be an approximation of the two-step method described above.
However, such an approximation may, nonetheless, provide enough accuracy to obtain good-quality clusters.
Lastly, it should be highlighted that the first step of the algorithm is completely amenable to be performed via parallel computation. 
The required time to carry out these evaluations can thus be greatly reduced.

The second step consists in minimizing a measure of clustering error, where the quantity in absolute value can be seen as a goodness of fit for cluster $C_k$. 
Minimizing this error pushes the cost function of the representative $r_k$ (represented by $V_{r_k,r_k}$) to fit as much as possible the average cost function of all the cluster's members (represented by $\tfrac{1}{\vert C_k\vert} \sum_{j\in C_k} V_{r_k,j}$). 
In other words, clusters are created such that, in each case, an elected representative can be identified that best fits the average ``view''.
It should be noted that this means that two scenarios with very different cost values may still be included in the same cluster if there exists a third scenario whose value provides a mid-ground between them. 
Consequently, it can be easily shown that the minimization problem \eqref{clustering_error} is equivalent to solving the mixed-integer program:
\begin{align}
\label{mip_obj}
\min & \quad \frac{1}{N} \sum_{i=1}^N t_i & \\
\label{mip_c1}
\text{s.t.} & \quad t_j \geq \sum_{i=1}^N x_{ij} V_{j,i} - \sum_{i=1}^N x_{ij} V_{j,j}, & \forall j\in\setz1N; \\
\label{mip_c2}
& \quad t_j \geq \sum_{i=1}^N x_{ij} V_{j,j} - \sum_{i=1}^N x_{ij} V_{j,i}, & \forall i\in\setz1N; \\
\label{mip_c3}
& \quad x_{ij} \leq u_j, \quad x_{jj} = u_j & \forall (i,j)\in\setz1N^2; \\
\label{mip_c4}
& \quad \sum_{j=1}^N x_{ij} = 1, \quad \sum_{j=1}^N u_j = K & \forall i \in \setz1N; \\
\label{mip_v}
& \quad x_{ij} \in \{0,1\}, \ u_{j} \in \{0,1\}, \ t_i \in [0,\infty), & \forall (i,j) \in \setz1N^2.
\end{align}

The binary variable $u_j$ determines if scenario $\xi_j$ is picked as a cluster representative, $\forall j \in\setz1N$, while the binary variable $x_{ij}$ establishes whether or not the scenario $\xi_i$ is in the cluster with representative $\xi_j$, $\forall (i,j)\in\setz1N^2$. 
The first two constraints linearize the absolute value of the objective function (i.e., the goodness of fit for each cluster). 
The other constraints guarantee that there are exactly $K$ clusters, that each one is non-empty and has a representative that belongs to it.

\subsection{Methodology Illustration}
\label{Section_method_illustration}

\begin{figure}
	\centering
	\includegraphics[scale=.5]{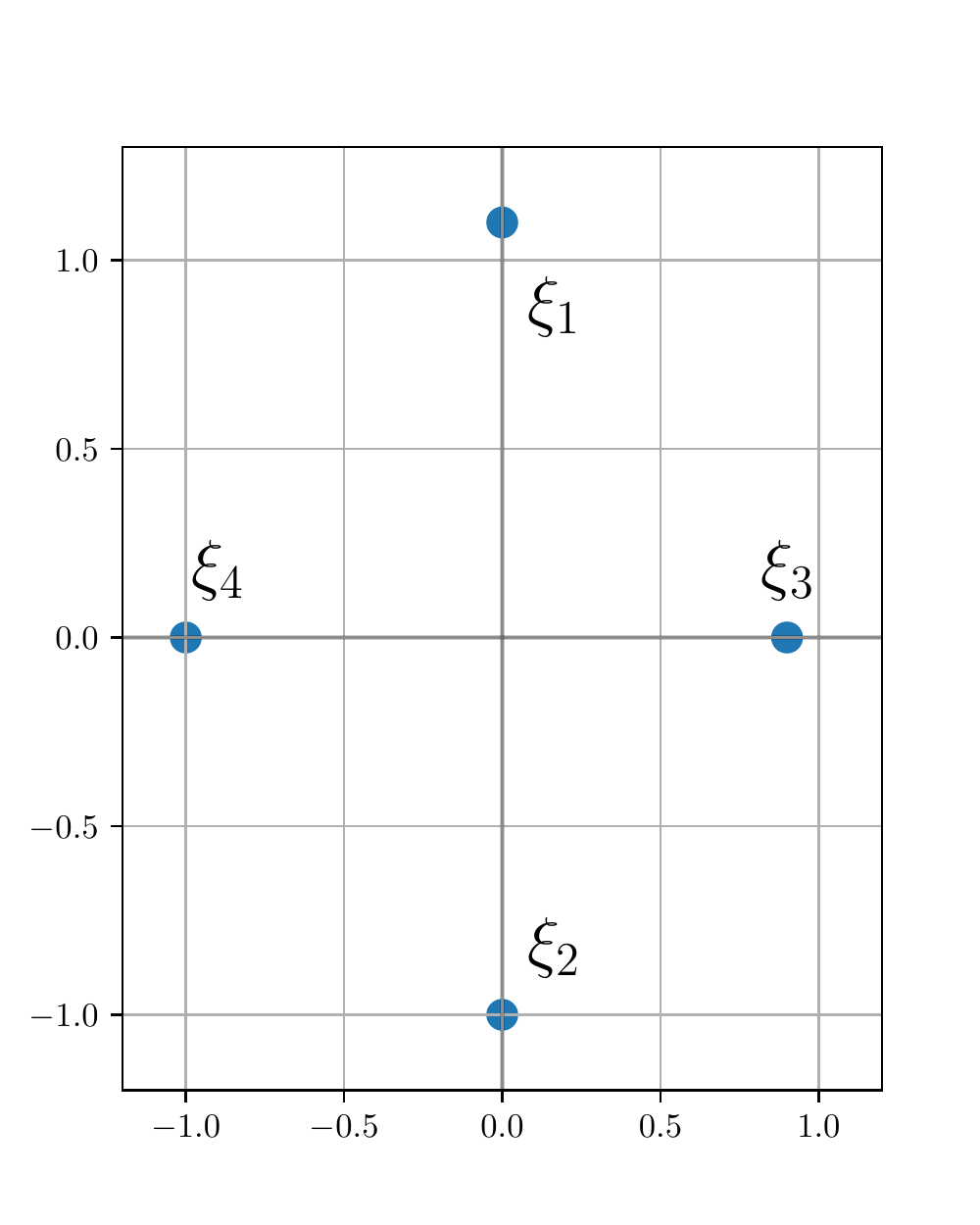}
	\caption{Location of the 4 scenarios of the stochastic problem described in \eqref{toy_problem_start}-\eqref{toy_problem_end}.}
	\label{fig_toy_problem_scenarios}
\end{figure}

\begin{figure}
	\centering
	\subfloat[$k$-means]{\includegraphics[scale=.5]{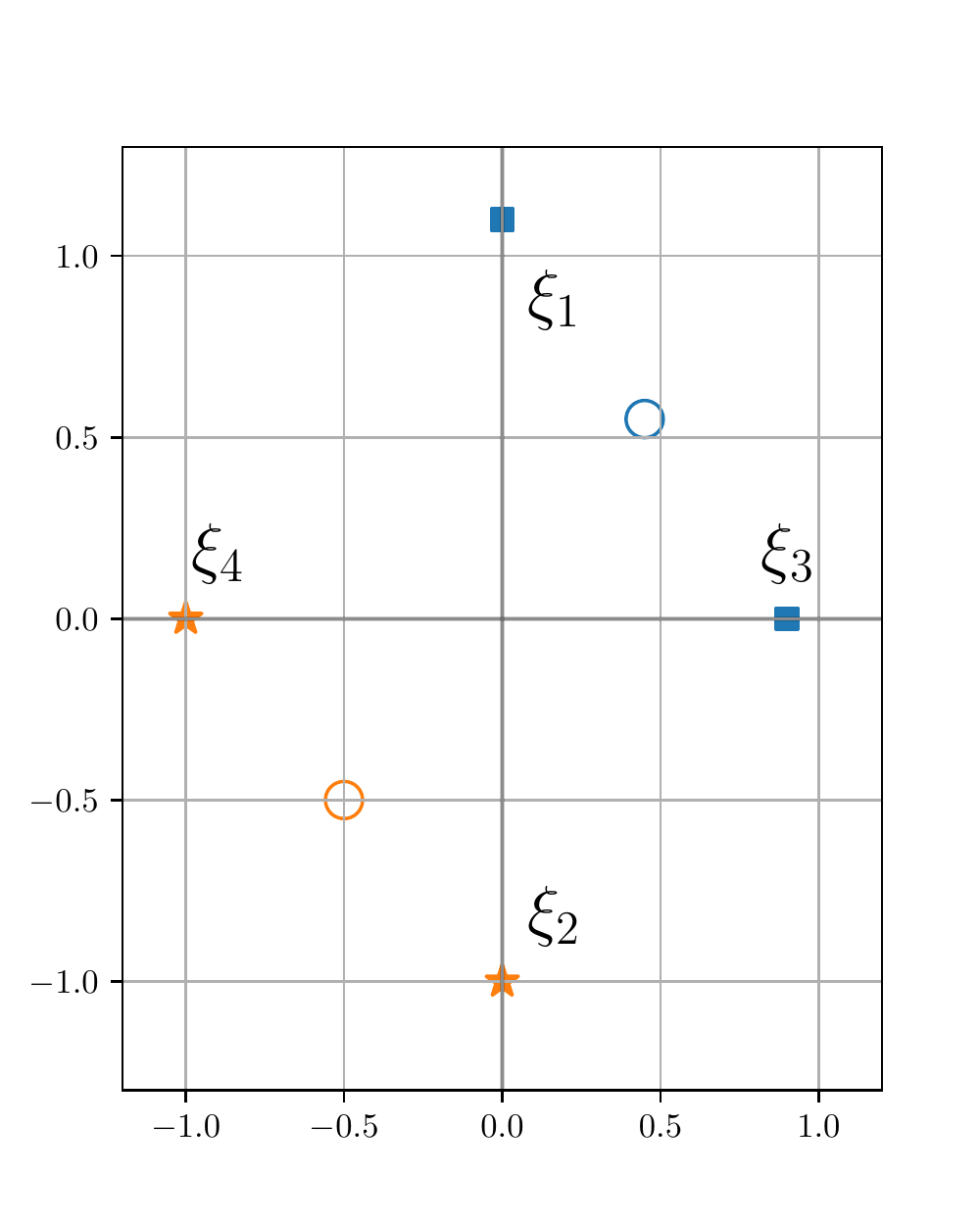}}
	\subfloat[$k$-medoids]{\includegraphics[scale=.5]{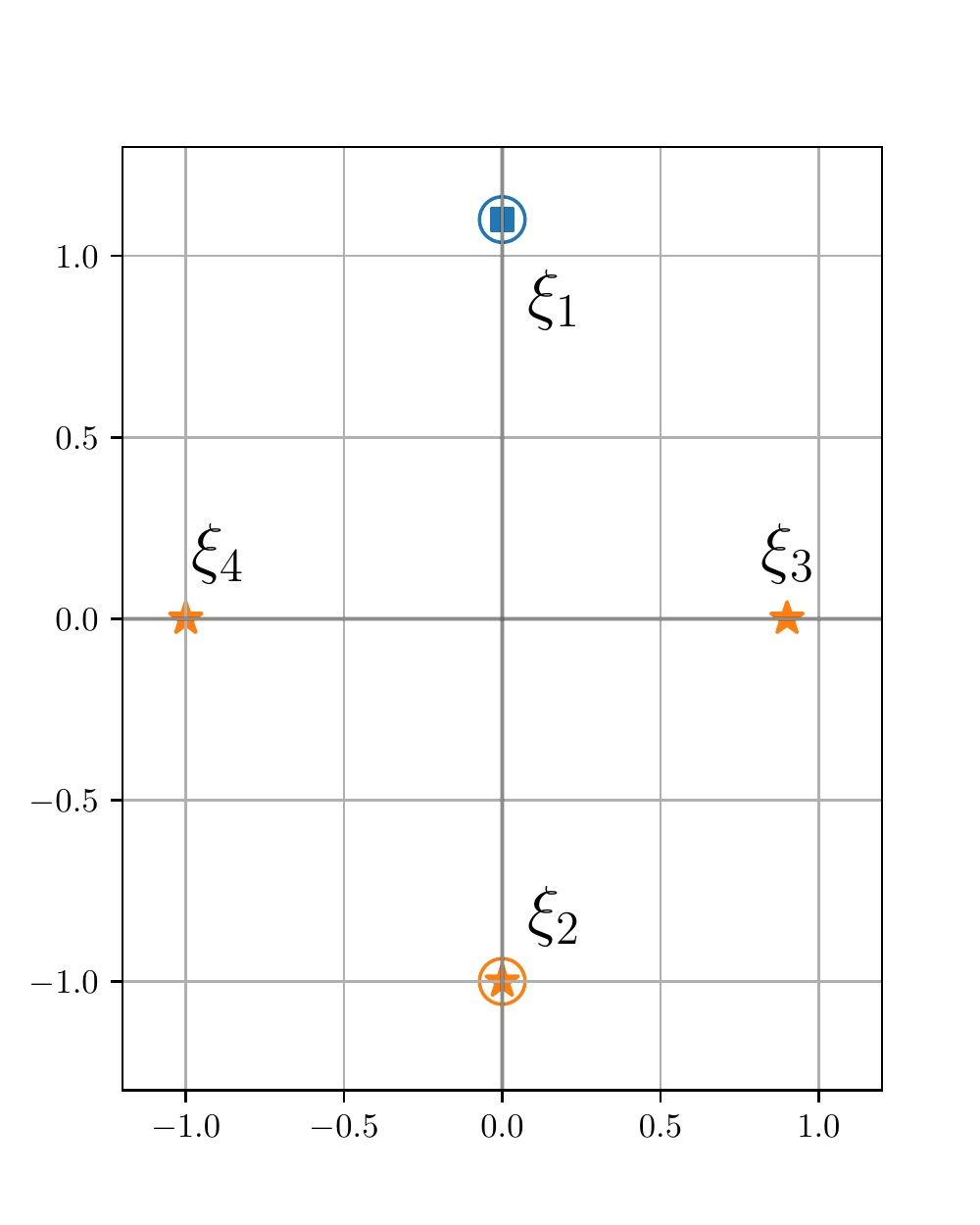}}
	\subfloat[CSSC]{\includegraphics[scale=.5]{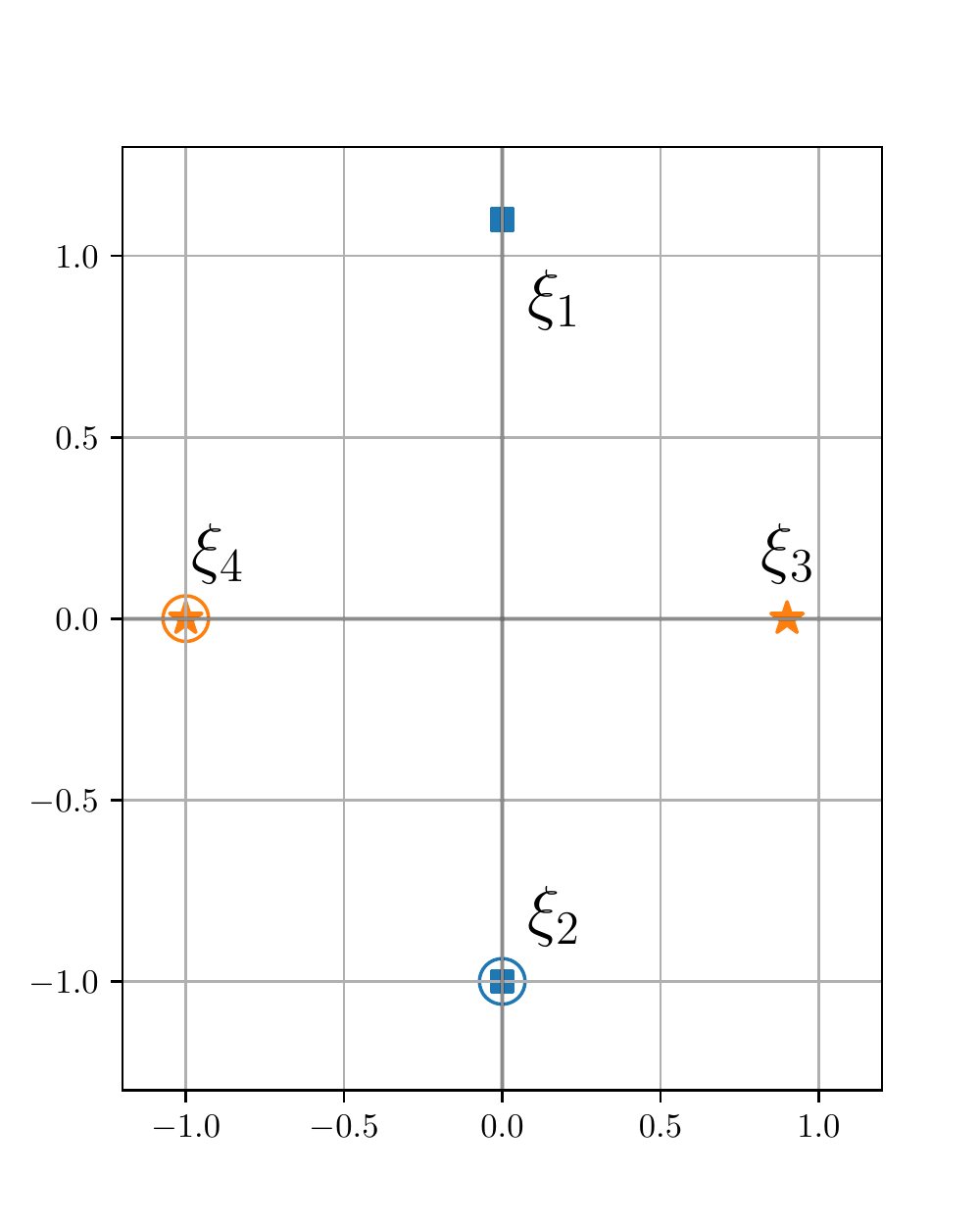}}
	\caption{Clustering performed by $k$-means, $k$-medoids and the CSSC algorithms for $K=2$. The representative of each cluster is displayed by a circle and the scenarios lying in the same cluster have the same color and marker.}
	\label{fig_toy_problem_clustering}
\end{figure}

To illustrate the CSSC algorithm and understand why it may provide an improvement over the purely distribution-based approaches, consider the following two-stage stochastic mixed-integer problem with four scenarios: 
\begin{align}
\label{toy_problem_start}
\underset{}{\min} & \ \frac{1}{4}\sum_{i=1}^4 \left(2t_1^i + 3t_2^i - y_1^i \xi_i^1 - y_2^i \xi_i^2\right) \\
\text{s.t.} & \ x + y_1^i \xi_i^1 - t_1^i \leq 0 & i\in\{1,2,3,4\}; \\
& \ x + y_1^i \xi_i^1 + t_1^i \geq 0 & i\in\{1,2,3,4\}; \\
& \ x - y_2^i \xi_i^2 - t_2^i \leq 0 & i\in\{1,2,3,4\}; \\
& \ x - y_2^i \xi_i^2 + t_2^i \geq 0 & i\in\{1,2,3,4\}; \\
\label{toy_problem_end}
& \ x \in \R, \ t_1^i \geq 0, \ t_2^i \geq 0, \ y_1^i\in\{-1,1\}, \ y_2^i\in\{-1,1\} & i\in\{1,2,3,4\}.
\end{align}
This problem has one decision variable at stage 0 ($x$), four decision variables at stage 1 ($t_1$, $t_2$, $y_1$, $y_2$), two random parameters ($\xi_i =(\xi_i^1, \xi_i^2)$), and four scenarios: $\xi_1 = (0, 0.9)$, $\xi_2 = (0, -1)$, $\xi_3 = (1.1, 0)$ and $\xi_4 = (-1, 0)$, which are displayed in Figure~\ref{fig_toy_problem_scenarios}. Its (unique) optimal solution is $x^* = 0$ with objective value $v^* = 1.475$.

Performing the $k$-means and $k$-medoids clustering methods on the four scenarios in Figure~\ref{fig_toy_problem_scenarios} for $K=2$ yields the clusters and representatives displayed in Figure~\ref{fig_toy_problem_clustering}(a)-(b). The $k$-means method finds two clusters of equal size and the representatives are the barycenters which lie in the middle of the straight line connecting each pair of points. Note that there exists another equivalent partition for $k$-means (which is not displayed) with the two clusters $\{\xi_1, \xi_4\}$ and $\{\xi_2,\xi_3\}$ (equivalence means that both have the same error, so the algorithm might output either one depending on the starting points). The partition generated by the $k$-medoids algorithm has a unique solution with scenarios $\{\xi_2, \xi_3, \xi_4\}$ together with representative $\xi_2$ and scenario $\{\xi_1\}$ alone, hence being its own representative. In the end, regardless of which algorithm is used, we see that the results follow the same pattern: scenarios close to each other in the space of random outcomes are clustered together. Of course, this is not a surprise, because these methods are designed for that purpose.

The CSSC algorithm performed on the same four scenarios yields a completely different result. Solving the problem on each scenario individually, as described in step 1 of the algorithm, produces the following opportunity cost matrix:
\begin{align}
\label{eq_toyV}
{\bm V} = 
\begin{bmatrix}
0.9 & 1.1 & 4.2 & 3.9 \\
1.4 & 1 & 4.3 & 4. \\
1.8 & 2 & 1.1 & 1. \\
1.8 & 2 & 1.1 & 1. 
\end{bmatrix}.
\end{align}
Recall that the $(i,j)$-element of $\bm{V}$ equals $F(x_i^*, \xi_j)$, i.e., the value of solution $x_i^*$ evaluated in scenario $\xi_j$, where $x_i^*$ is the optimal solution of the deterministic problem built from the individual scenario $\xi_i$.

Applying now step 2 of the algorithm, i.e., solving the MIP formulation \eqref{mip_obj}-\eqref{mip_v} for the matrix \eqref{eq_toyV} and $K=2$, produces the partition and representatives displayed in Figure~\ref{fig_toy_problem_clustering}(c). We see that it looks somehow counter-intuitive (at least when we look at it in the space of random outcomes), since it pairs together scenarios that are the furthest away from each other. There is, however, a good reason for that, which lies in the structure of the stochastic problem, and which is reflected in the opportunity cost matrix \eqref{eq_toyV}: the first two columns of the matrix have similar values, as well as the last two columns. This means that the cost function $F(x, \xi_i)$ satisfies:
\begin{equation}
\label{equa_proximity_of_costs}
F(x, \xi_1) \simeq F(x, \xi_2) \quad \text{ and } \quad F(x, \xi_3) \simeq F(x, \xi_4), \quad \forall x\in \{x_1^*, x_2^*, x_3^*, x_4^*\},
\end{equation}
where $\{x_1^*, x_2^*, x_3^*, x_4^*\}$ is the set of first-stage solutions computed in \eqref{eq_one_scenario_problem} for each individual scenarios. This proximity of scenarios in the space of cost values suggests pairing together scenarios $\{\xi_1, \xi_2\}$ and $\{\xi_3,\xi_4\}$, even though they appear to be the furthest away in the space of random outcomes. The basis for pairing them together, which was already explained in section~\ref{motivation_cost_driven}, is that two scenarios $\xi_i$ and $\xi_j$ can be merged without error if $F(x,\xi_i) = F(x,\xi_j)$ for all feasible $x\in X$. Although in \eqref{equa_proximity_of_costs} the equality does not hold exactly but only approximately over a subset of feasible solutions, which means that some errors might still be introduced by doing such scenario reduction, we can expect a smaller error to be introduced by following the guidelines provided by the opportunity cost matrix.

This expectation turns out to be true in the specific problem considered here, since the approximate problem solved using the $k$-means and $k$-medoids scenarios provide the solutions $x=0.5$ of true value 2.475 ($62\%$ away from optimality) and $x=1$ of true value 2.575 ($69\%$ away from optimality), respectively, while the CSSC scenarios do provide the actual optimal solution $x^* = 0$. An exact scenario reduction can thus be achieved for this problem, not by chance, but because of an underlying property hidden in the problem's objective function and constraints. Indeed, one could see by computing exactly or numerically the cost function $F(x, \xi)$ for any $\xi=(\xi^1,\xi^2)\in\R^2$ that it approximately satisfies
\begin{equation}
F(x,\xi) \simeq 2\left\vert x + \vert \xi^1\vert \right\vert + 3\left\vert x - \vert \xi^2 \vert\right \vert - \vert \xi^1 \vert - \vert \xi^2 \vert.
\end{equation}
Thus, $F$ depends on $(\xi^1, \xi^2)$ only through the absolute values $(\vert \xi^1 \vert, \vert \xi^2 \vert)$. This explains why scenarios at the opposite of each other in Figure~\ref{fig_toy_problem_scenarios} provide similar cost values, and hence can be clustered together. This property is of course invisible for the clustering methods that only work in the space of distributions. 

\section{Numerical Experiments}
\label{sec:numerical}

In this section we apply the CSSC algorithm on two stochastic mixed-integer programs:  1) network design and  2) facility location problems. 
The goals of the conducted experiments are twofold.
Using the stochastic network design programs, we study the problem of scenario reduction in general (i.e., for varying values of $K$).
While in the case of the stochastic facility location programs, we focus on how good a solution with a \emph{single} scenario can be.
Considering a unique scenario means turning the stochastic problem into a deterministic one, which is the path followed by many decision-makers when the problem becomes too hard to be solved as a stochastic one. 
The facility location problem, which has a number of second stage integer decisions typically several order of magnitude larger than the one in the first stage, makes a good candidate for such extreme scenario reduction. 
Our computational experiments will show that, in this situation, the CSSC algorithm also provides a valuable deterministic scenario estimate that outperforms estimates based on the scenarios alone.

All experiments in this section are run on Python 3.7 and the mixed-integer problems are solved using CPLEX Optimization Studio version 12.9 via its Python API. The network design problem is solved on a Windows 10 computer with 4 cores @2.11GHz 16GB RAM, and the facility location problem is solved on a Windows 10 computer with 6 cores @3.2GHz 32GB RAM.

We compare the CSSC algorithm to several clustering methods, which are all available open-source. The $k$-means algorithm is implemented in the \emph{scikit-learn} library (see \cite{scikit_learn}). The $k$-medoids and $k$-medians implementations are those of the \emph{pyclustering} library (see \cite{novikov2019}). For the $k$-modes algorithm, we refer to \cite{devos2015}.
 
\subsection{Network Design Problem}
\label{network_design}

We first consider a stochastic network design problem, where a number of commodities are to be transported across a network that has to be designed before the demand of these commodities is known. 
The problem is formulated as a two-stage model. 
In the first stage, the decisions about the structure of the network are made (i.e., which arcs are opened and which are closed)
Then, the stochastic demands are revealed and, based on them, the second stage decisions of commodity flows in the network are made.
Considering its wide applicability, ranging from supply chain management, e.g., \cite{Santoso:EJOOR:2005}, to network telecommunication planning, e.g., \cite{Riis:IJOC:2002},  there has been a steady stream of research dedicated to the development of solution methods for this problem, see \cite{Crainic:COR:2014} and \cite{Rahmaniani:SJO:2018} for examples of general heuristics and an exact method, respectively.
However, considering their combinatorial nature, these problems remain extremely challenging to solve, especially when considering large networks and a high number of scenarios used to define the formulation.

Our setting is therefore a complete directed graph $G=(V,A)$, across which commodities from a set $C$ must be transported. 
Each arc $a\in A$, pointing from $a(0)\in V$ to $a(1)\in V$, has a total capacity $u_a$ if it has been opened for a fixed cost $c_a$. 
The cost of transporting one unit of commodity $c\in C$ across arc $a\in A$ is denoted by $q_{ac}$. 
We denote the demand for commodity $c\in C$ at vertex $v\in V$ under scenario $i$ by $d_{v,c}^i$. 

\begin{figure}[h]
	\centering
	\subfloat[1st-stage]{\includegraphics[trim={2cm 1cm 1.5cm 0cm}, clip, scale=0.46]{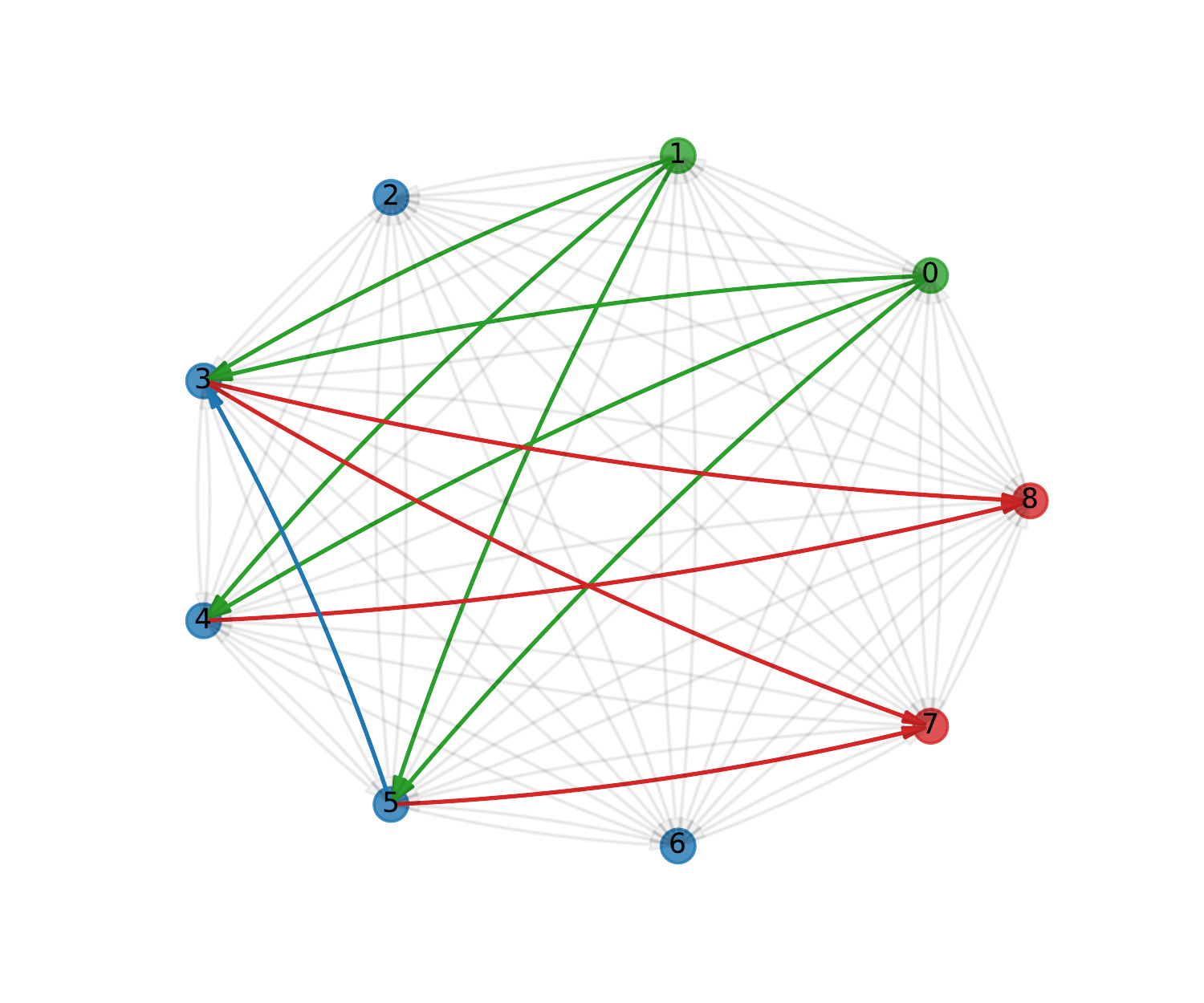}}\\%
	\subfloat[2nd-stage, scenario\#1]{\includegraphics[trim={2cm 1cm 1.5cm 0.2cm}, clip, scale=0.46]{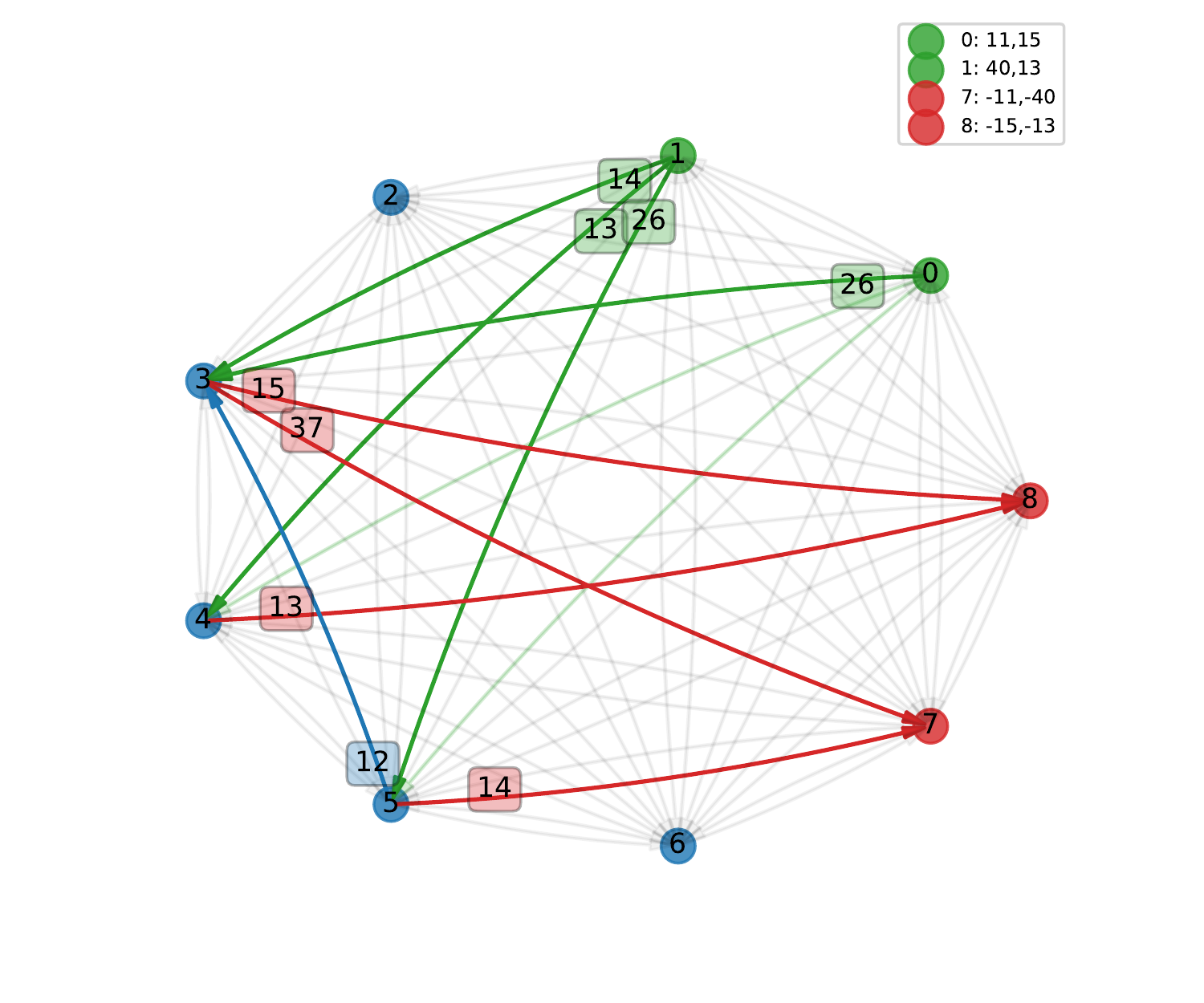}}%
	\qquad
	\subfloat[2nd-stage, scenario\#2]{\includegraphics[trim={2cm 1cm 1.5cm 0.2cm}, clip, scale=0.46]{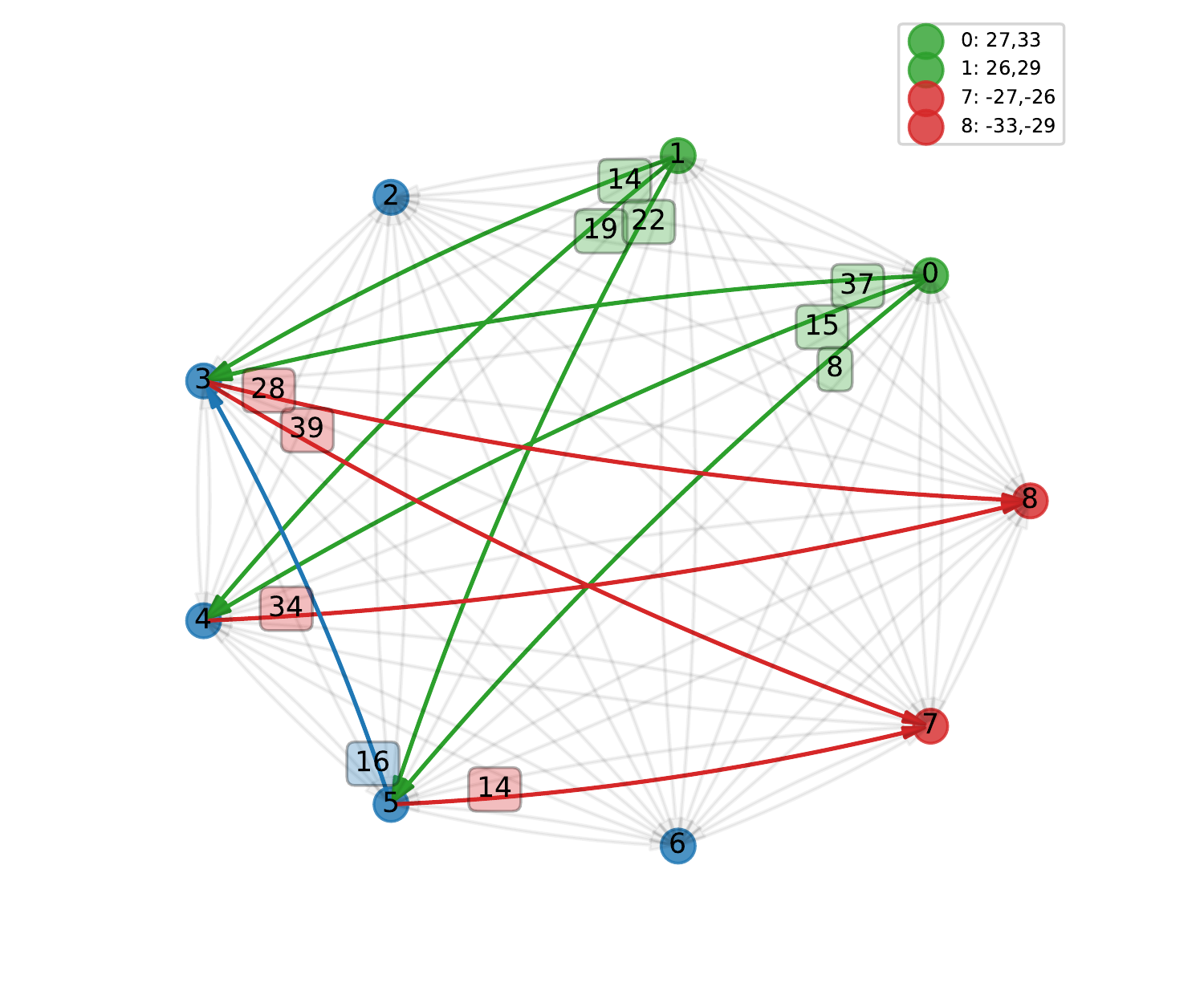}}%
	\caption{First and second stage decisions for a network design problem with 4 commodities, 9 vertices, and 68 arcs. The commodities are transported from the green vertices (origins) to the red ones (destinations), passing through the blue ones (intermediates). Each pair (origin, destination) is a commodity whose random demand is displayed at the top right corner of the network. Arcs are colorized to increase visibility: they are green if they leave an origin, red if they enter a destination, and blue if they connect two intermediates vertices.}%
	\label{fig:NDP}%
\end{figure}

The first stage decision of opening the arc $a\in A$ is denoted by $x_a$, which can take two values: 0 if the arc is closed and 1 if it is open. The second stage decision $y^i_{ac}$ is the number of units of commodity $c\in C$ transported across arc $a\in A$ under scenario $i$. 
This problem is solved using the following two-stage stochastic formulation:
\begin{align}
\min & \ \sum_{a\in A} c_{a} x_{a} + \frac{1}{N}\sum_{i=1}^N \sum_{c\in C} \sum_{a\in A} q_{ac} y_{ac}^{i} \label{SND-obj}\\
\text{s.t.} & \ \sum_{\substack{a\in A \\ a(0)=v}} y_{ac}^{i} - \sum_{\substack{a\in A \\ a(1)=v}} y_{ac}^{i} = d_{v,c}^{i} & \forall (v, c, i) \in V \times C \times \setz1N  \label{SND-flow}\\
& \ \sum_{c\in C} y_{ac}^{i} \leq u_{a} x_{a} & \forall (a, i) \in A\times\setz1N  \label{SND-capacity}\\
& \ x_{a} \in \{0,1\}, \ y_{ac}^{i} \in [0,\infty) & \forall (a,c,i) \in A\times C \times \setz1N.  \label{SND-integrality-non-negativity}
\end{align}
The goal of the problem is to satisfy all demands while minimizing the expected total cost incurred, objective function (\ref{SND-obj}), which includes the fixed cost of opening the arcs plus variable cost of transporting the commodities across the network.
The constraints (\ref{SND-flow}) are the flow conservation requirements imposed for each vertex of the network, for each commodity and for all scenarios.
The capacities associated to the arcs, which limit the total amount of commodity flow that can be transported through them under each scenario, are enforced through the inclusion of constraints (\ref{SND-capacity}).
Finally, constraints (\ref{SND-integrality-non-negativity}) impose the necessary integrality and non-negativity requirements on the decision variables.
We refer to Figure~\ref{fig:NDP} for an illustration of first and second stage decisions for a network with 4 commodities, 9 vertices and 68 arcs.

To test the CSSC algorithm, we generate randomly two batches of instances. The first one contains 25 instances of a network with 49 commodities, 763 arcs, and 25 demand scenarios drawn from a discrete uniform distribution in $U(0,10)$. The second batch contains 100 instances of a network with 9 commodities, 132 arcs and 75 scenarios sampled from a log-normal distribution with mean 1 and variance 1 (which is then rounded-off to the nearest integer). Each instance is generated by sampling the (deterministic) parameters $q_{ac}$, $c_a$, $u_a$ independently and uniformly from fixed sub-intervals of the integers ($q_{ac}\sim U[5,10]$, $c_a\sim U[3,10]$, $u_a\sim U[10,40]$).

Four comparison methods are considered for this problem: $k$-means, $k$-medoids, $k$-medians and Monte Carlo sampling. For each one, the approximate problem \eqref{eq_approximate_problem} is solved for a reduced number of scenarios $K=2,\dots,10$. The resulting optimal solution is then evaluated in the original problem to assess its true value and the gap with the optimal value of the original problem (as defined in \eqref{implementation_error}).
The implementation error thus serves as the quality measure to compare the methods. 
Note that this error is always expressed in percentage of the original optimal value to allow comparison between different instances. 

\begin{figure}[h]
	\centering
	\includegraphics[trim={0 0 0 0}, clip, scale=0.75]{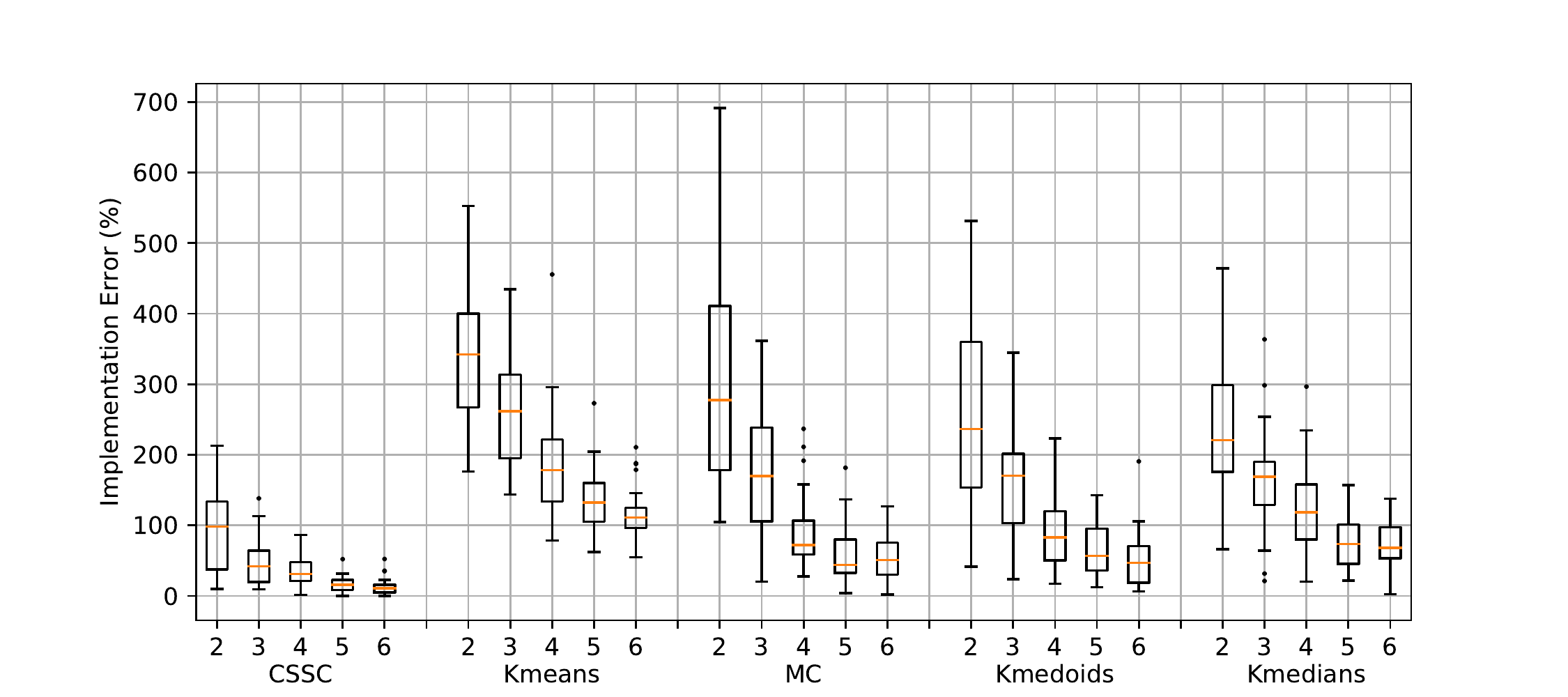}\\
	\includegraphics[trim={0 0 0 0}, clip, scale=0.75]{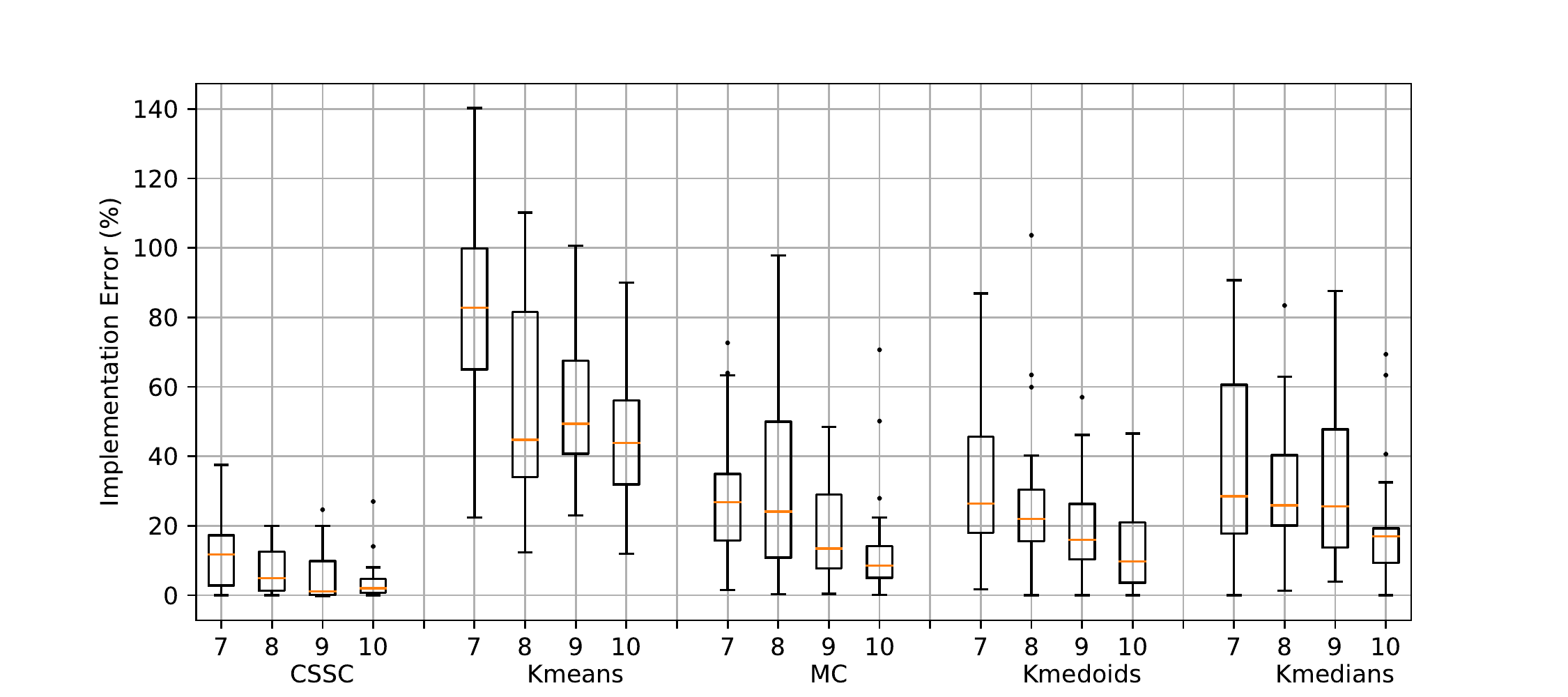}
	\caption{Implementation error (in percentage) for the network design problem with 25 scenarios. The comparison is on the range $K=2,\dots,6$ (top) and $K=7,\dots,10$ (bottom).}%
	\label{fig:NDP-25scen}%
\end{figure}

\begin{figure}[h]
	\centering
	\includegraphics[trim={0 0 0 0}, clip, scale=0.75]{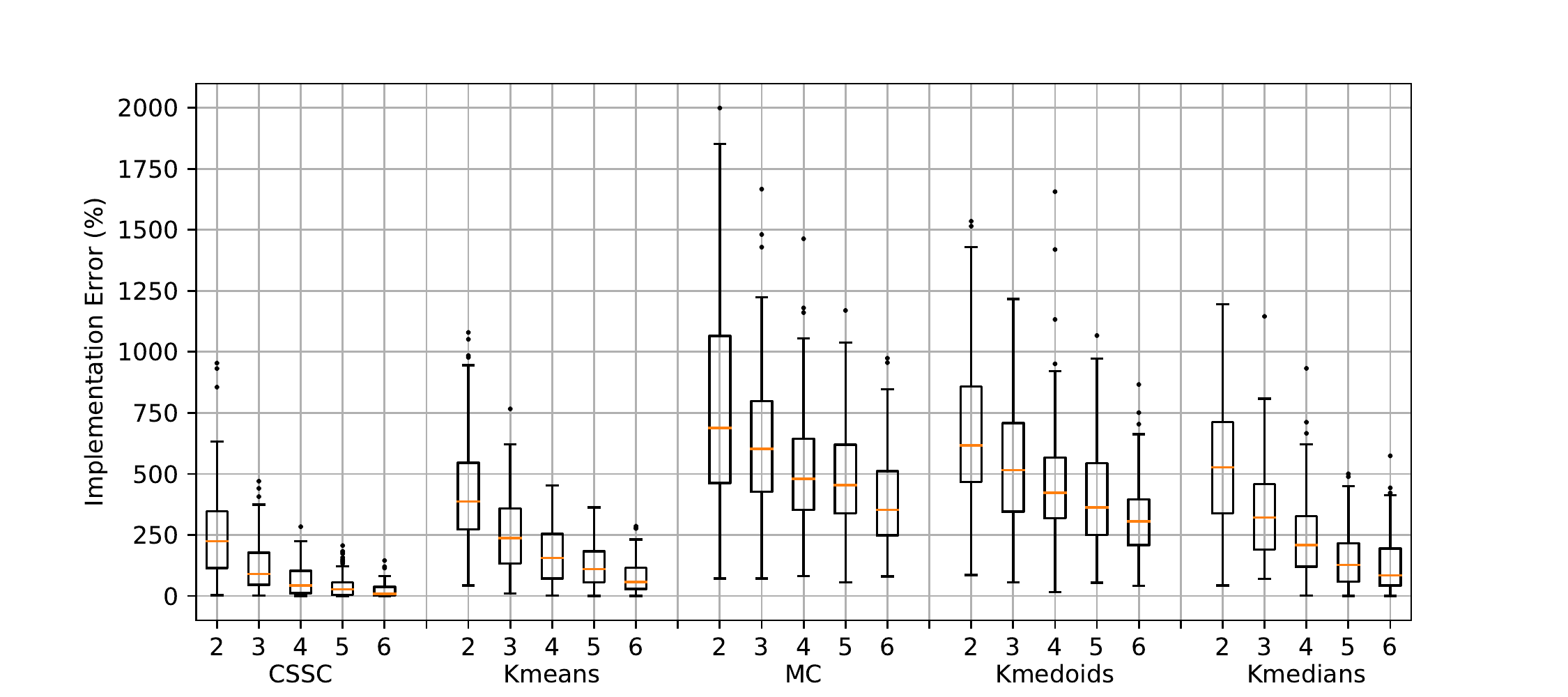}\\
	\includegraphics[trim={0 0 0 0}, clip, scale=0.75]{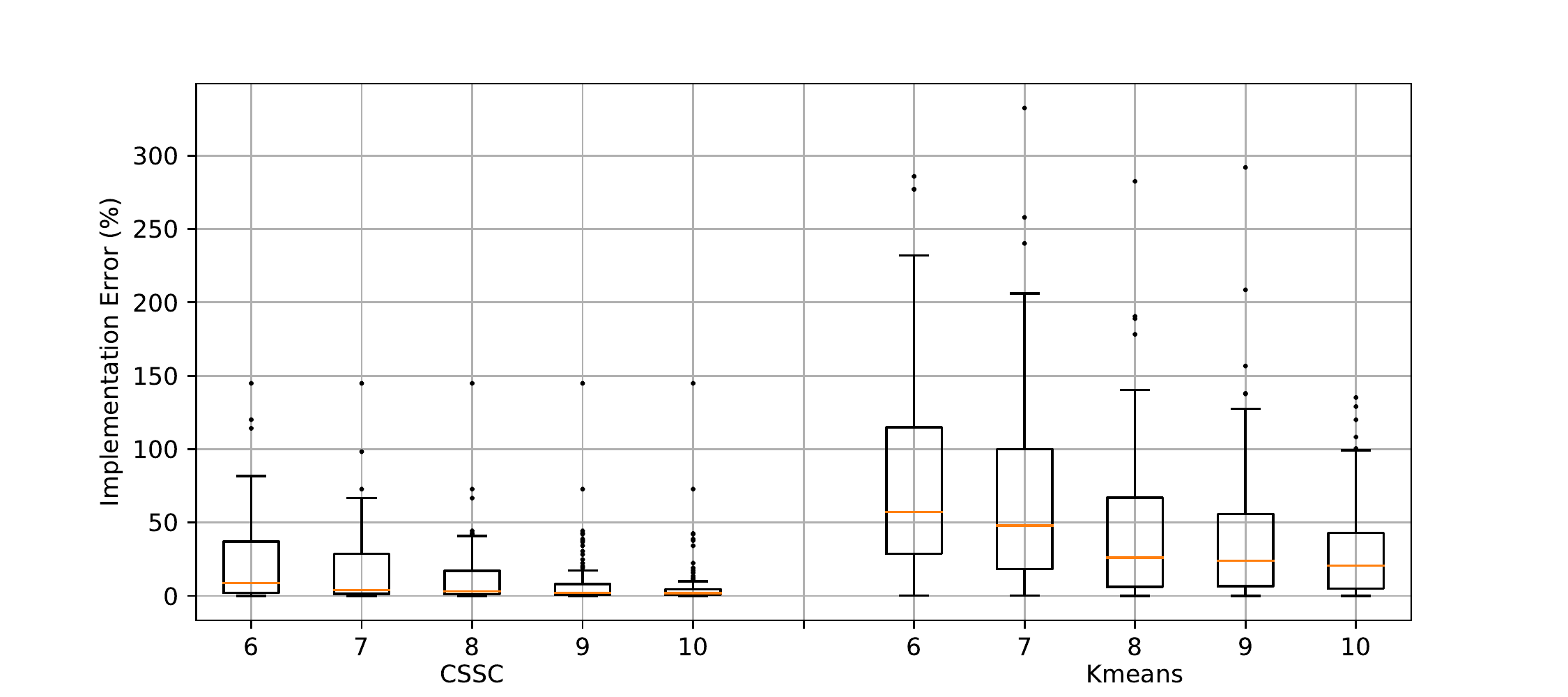}
	\caption{Implementation error for the network design problem with 75 scenarios. The comparison is over all methods on the range $K=2,\dots,6$ (top) and only between CSSC and $k$-means on the range $K=6,\dots,10$ (bottom) because the other methods have much larger errors on that range.}%
	\label{fig:NDP-75scen}%
\end{figure}

The results are displayed as box plots in Figure~\ref{fig:NDP-25scen} for the instances with 25 scenarios and in Figure~\ref{fig:NDP-75scen} for the instances with 75 scenarios. At the top of each figure, all methods are compared and at the bottom only the best ones are displayed to ease the comparison, as some methods still yield large implementation errors that would not fit within the range of values displayed. The convention used for the box plots is the following: the line of the box represents the median, the box itself spans over the interquartile range of the data (i.e., the first and third quartile), the whiskers that extend from the box correspond to 1.5 times the interquartile range, and everything beyond the whiskers is plotted as individual points.

It appears clearly from these figures that the CSSC algorithm outperforms all the other comparison methods on every single clustering size for $K=2$ to $K=10$. For $K=10$, the median error becomes extremely close to zero and is several times smaller than even the best of the other methods. To quantify more precisely how many instances are solved with little error, we show in Table~\ref{Table:NDP-25scen} and \ref{Table:NDP-75scen} the percentage of instances solved with less than 10\% and 2\% error. For the instances with 25 scenarios and uniform stochastic demands, we see that more than 90\% of instances are solved with less than 10\% error by CSSC and about 50\% with less than 2\% error. The second best method, Monte Carlo, only achieves 60\% and 16\%, respectively, on those instances. The situation is similar for the instances with 75 scenarios and lognormal stochastic demands: 86\% of instances are solved with less than 10\% error and about 50\% with less than 2\%, against 36\% and 16\%, respectively, for the second best method ($k$-means in that case).

\begin{figure}[h]
	\centering
	\includegraphics[trim={0 0 0 0}, clip, scale=0.75]{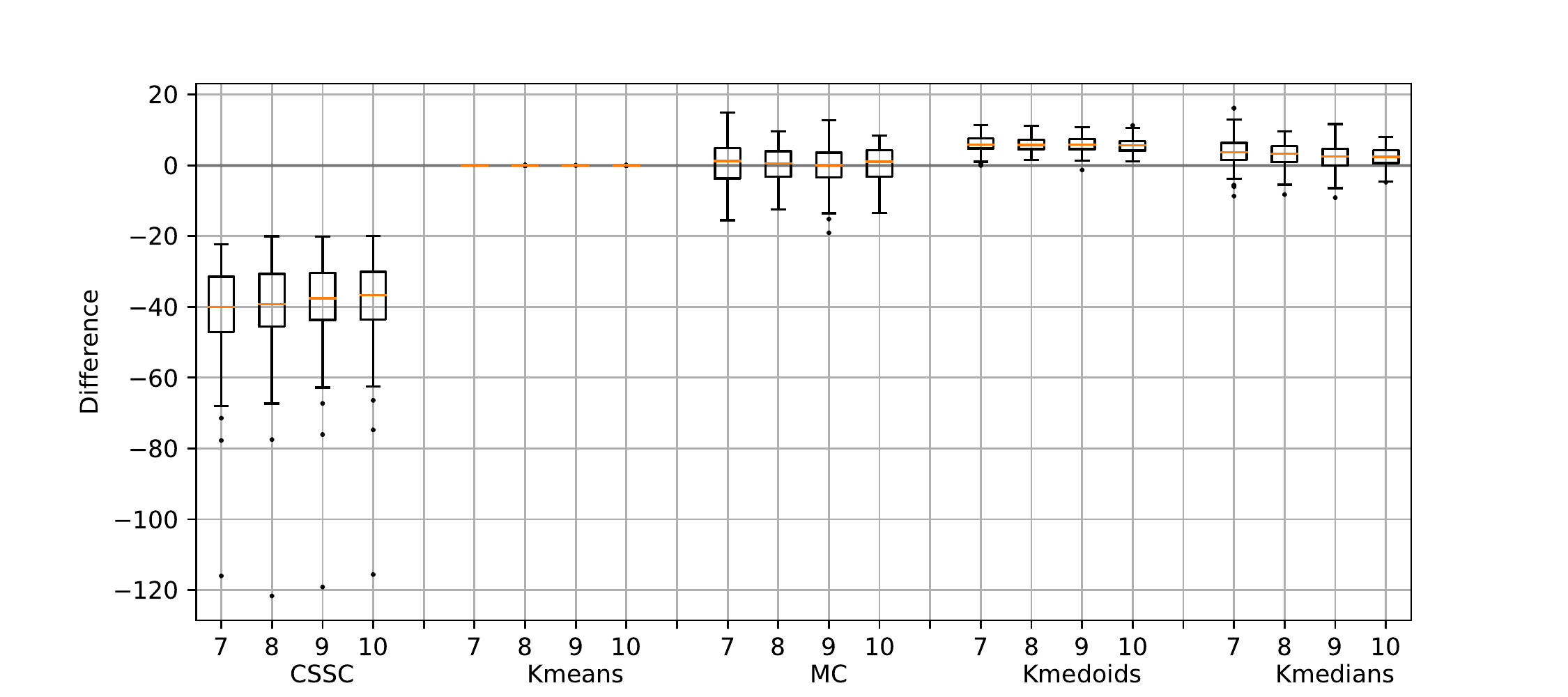}\\
	\includegraphics[trim={0 0 0 0}, clip, scale=0.75]{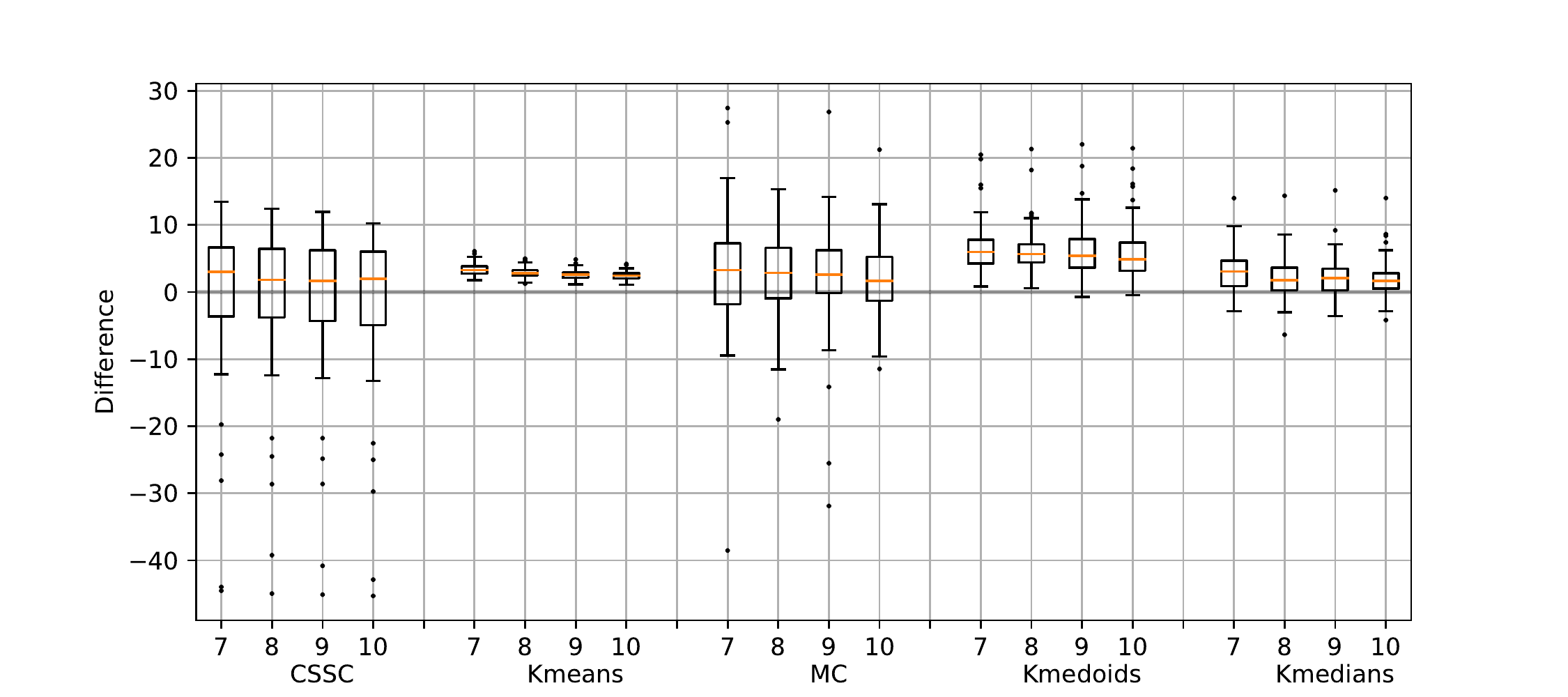}
	\caption{Difference between the statistic of the original set of scenarios and the reduced one for the mean (top) and the standard deviation (bottom) of the total demand in the network for the instances with 75 scenarios.}%
	\label{fig:NDP-75scen-mean-std}%
\end{figure}

\begin{figure}[h]
	\centering
	\includegraphics[trim={0 0 0 0}, clip, scale=0.75]{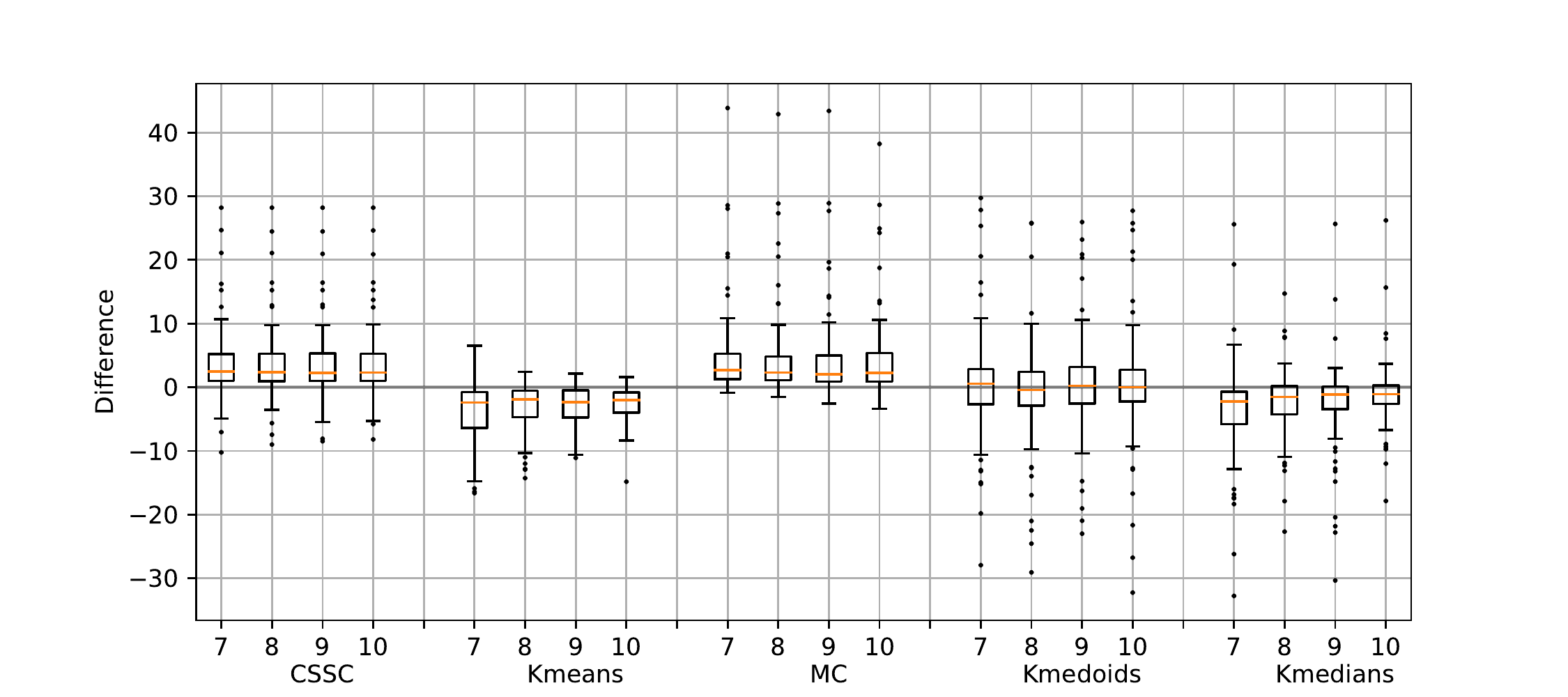}\\
	\includegraphics[trim={0 0 0 0}, clip, scale=0.75]{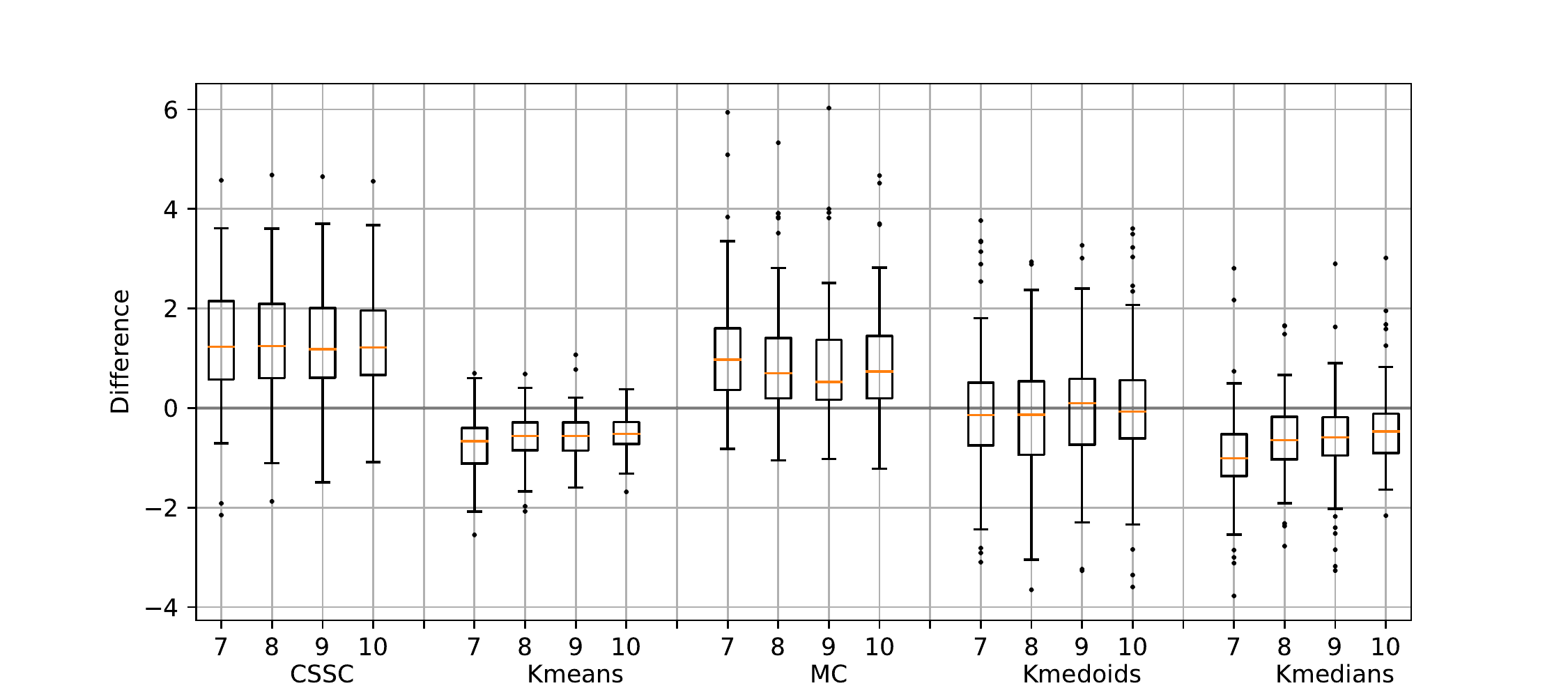}
	\caption{Difference between the statistic of the original set of scenarios and the reduced one for the kurtosis (top) and the skewness (bottom) of the total demand in the network for the instances with 75 scenarios.}%
	\label{fig:NDP-75scen-kurt-skew}%
\end{figure}

Now that we have seen that the CSSC algorithm is best at selecting a reduced set of scenarios to solve the network design problem, it is interesting to study the properties of those scenarios selected by CSSC but not by the other methods. To this end, we focus on the statistics of the total demand in the network, i.e., $\sum_{v,c} d_{v,c}^{i}$. Surprisingly, we observe that the scenarios of the CSSC algorithm are the worst at matching the mean of the total demand, as displayed at the top of Figure~\ref{fig:NDP-75scen-mean-std}. Indeed, while $k$-means matches it perfectly (which explains why no box is apparent), and the other methods are close to a perfect match within some statistical errors, the CSSC algorithm completely overestimates the mean demand in the network for every single instance. The other moments (variance, skewness, kurtosis) are also not properly matched by the scenarios of the CSSC algorithm, although the differences in those cases are less dramatic (see Figure~\ref{fig:NDP-75scen-mean-std} and \ref{fig:NDP-75scen-kurt-skew}). We note that this goes against the guidelines generally followed in stochastic programming that scenarios should match the moments of the distribution (see, e.g., \cite{chopra2013effect}). Of course, this observation is limited to the specific problem considered here, and is not an evidence that the moments of the distribution should not be matched, but rather that matching them may not always be necessary to achieve good approximation results.

As for the computational times to run the CSSC algorithm, the first step of the algorithm (computation of the opportunity cost matrix) takes between 18 sec and 45 sec (mean: 23 sec) for the instances with 25 scenarios and between 4 sec and 11 sec (mean: 8.2 sec) for the instances with 75 scenarios. The latter takes less time despite the higher number of scenarios because the network is smaller in that case. Note that we have not used any parallelization scheme to speed up the computation. We also recall that the opportunity cost matrix is computed only once per instance, and not for each value of $K$ used in the clustering. As for the second step (solution of the MIP partitioning problem), it takes always less than 1 sec for 25 scenarios and between 2 sec and 4 sec (mean: 2.9 sec) for 75 scenarios. Finally, these computation times are fairly constant over the range of tested $K$. 

\begin{table}[h]
	\centering
	\begin{tabular}{|c|ccccc|ccccc|}
		\hline
		& \multicolumn{5}{|c|}{$\mathrm{Error} \leq 10\%$} & \multicolumn{5}{|c|}{$\mathrm{Error} \leq 2\%$}\\
		\hline
$K$ & CSSC & Kmeans & MC & Kmedoids & Kmedians & CSSC & Kmeans & MC & Kmedoids & Kmedians \\
\hline
1 & 	 0 & 	 0 & 	 0 & 	 0 & 	 0 & 	 0 & 	 0 & 	 0 & 	 0 & 	 0  \\
2 & 	 4 & 	 0 & 	 0 & 	 0 & 	 0 & 	 0 & 	 0 & 	 0 & 	 0 & 	 0  \\
3 & 	 4 & 	 0 & 	 0 & 	 0 & 	 0 & 	 0 & 	 0 & 	 0 & 	 0 & 	 0  \\
4 & 	 4 & 	 0 & 	 0 & 	 0 & 	 0 & 	 4 & 	 0 & 	 0 & 	 0 & 	 0  \\
5 & 	 28 & 	 0 & 	 8 & 	 0 & 	 0 & 	 8 & 	 0 & 	 0 & 	 0 & 	 0  \\
6 & 	 40 & 	 0 & 	 8 & 	 4 & 	 4 & 	 16 & 	 0 & 	 0 & 	 0 & 	 0  \\
7 & 	 48 & 	 0 & 	 16 & 	 12 & 	 12 & 	 12 & 	 0 & 	 8 & 	 4 & 	 4  \\
8 & 	 60 & 	 0 & 	 20 & 	 20 & 	 12 & 	 32 & 	 0 & 	 4 & 	 4 & 	 4  \\
9 & 	 76 & 	 0 & 	 28 & 	 24 & 	 16 & 	 56 & 	 0 & 	 4 & 	 4 & 	 0  \\
10 & 	 92 & 	 0 & 	 60 & 	 52 & 	 28 & 	 48 & 	 0 & 	 16 & 	 12 & 	 8  \\
		\hline
	\end{tabular}
	\caption{Percentage of instances solved with less than $10\%$ and $2\%$ error for the network design problem with 25 scenarios.}
	\label{Table:NDP-25scen}
\end{table}

\begin{table}[h]
	\centering
\begin{tabular}{|c|ccccc|ccccc|}
	\hline
	& \multicolumn{5}{|c|}{$\mathrm{Error} \leq 10\%$} & \multicolumn{5}{|c|}{$\mathrm{Error} \leq 2\%$}\\
	\hline
$K$ & CSSC & Kmeans & MC & Kmedoids & Kmedians & CSSC & Kmeans & MC & Kmedoids & Kmedians \\
\hline
1 & 	 0 & 	 0 & 	 0 & 	 0 & 	 0 & 	 0 & 	 0 & 	 0 & 	 0 & 	 0 \\
2 & 	 3 & 	 0 & 	 0 & 	 0 & 	 0 & 	 0 & 	 0 & 	 0 & 	 0 & 	 0 \\
3 & 	 8 & 	 0 & 	 0 & 	 0 & 	 0 & 	 3 & 	 0 & 	 0 & 	 0 & 	 0 \\
4 & 	 22 & 	 2 & 	 0 & 	 0 & 	 1 & 	 9 & 	 1 & 	 0 & 	 0 & 	 1 \\
5 & 	 29 & 	 6 & 	 0 & 	 0 & 	 5 & 	 14 & 	 5 & 	 0 & 	 0 & 	 2 \\
6 & 	 51 & 	 11 & 	 0 & 	 0 & 	 10 & 	 25 & 	 4 & 	 0 & 	 0 & 	 7 \\
7 & 	 59 & 	 17 & 	 1 & 	 0 & 	 13 & 	 32 & 	 7 & 	 1 & 	 0 & 	 6 \\
8 & 	 66 & 	 30 & 	 0 & 	 0 & 	 17 & 	 39 & 	 13 & 	 0 & 	 0 & 	 9 \\
9 & 	 78 & 	 32 & 	 0 & 	 1 & 	 21 & 	 47 & 	 10 & 	 0 & 	 1 & 	 11 \\
10 & 	 86 & 	 36 & 	 1 & 	 2 & 	 29 & 	 51 & 	 16 & 	 1 & 	 0 & 	 16 \\
		\hline
	\end{tabular}
	\caption{Percentage of instances solved with less than $10\%$ and $2\%$ error for the network design problem with 75 scenarios.}
	\label{Table:NDP-75scen}
\end{table}

\subsection{Facility Location Problem}
\label{facility_location}

This two-stage problem consists in locating facilities over a set of predefined positions and assigning them to randomly present customers. The decision of positioning (or not) a facility at some location is made at the first stage, before knowing how many customers will be present and where. Then, the scenario of customers presence is revealed and, based on it, the facilities are matched to the customers in the most cost-efficient way. 

The facility location problem has been well-studied in the literature in its deterministic form (e.g., in \cite{guha1999greedy, li20131}, among many others) and to a much smaller extent in a stochastic version with random customers demand (e.g., \cite{bieniek2015note, baron2008facility}). The problem that we consider here does not feature random demand but \emph{random presence} of customers. The following formulation is inspired by that of \cite{ntaimo2005million} where the problem is referred to as the server location problem.

\begin{figure}[h]
	\centering
	\subfloat[1st-stage]{\includegraphics[trim={1cm 1cm 1cm 1cm}, clip, scale=0.5]{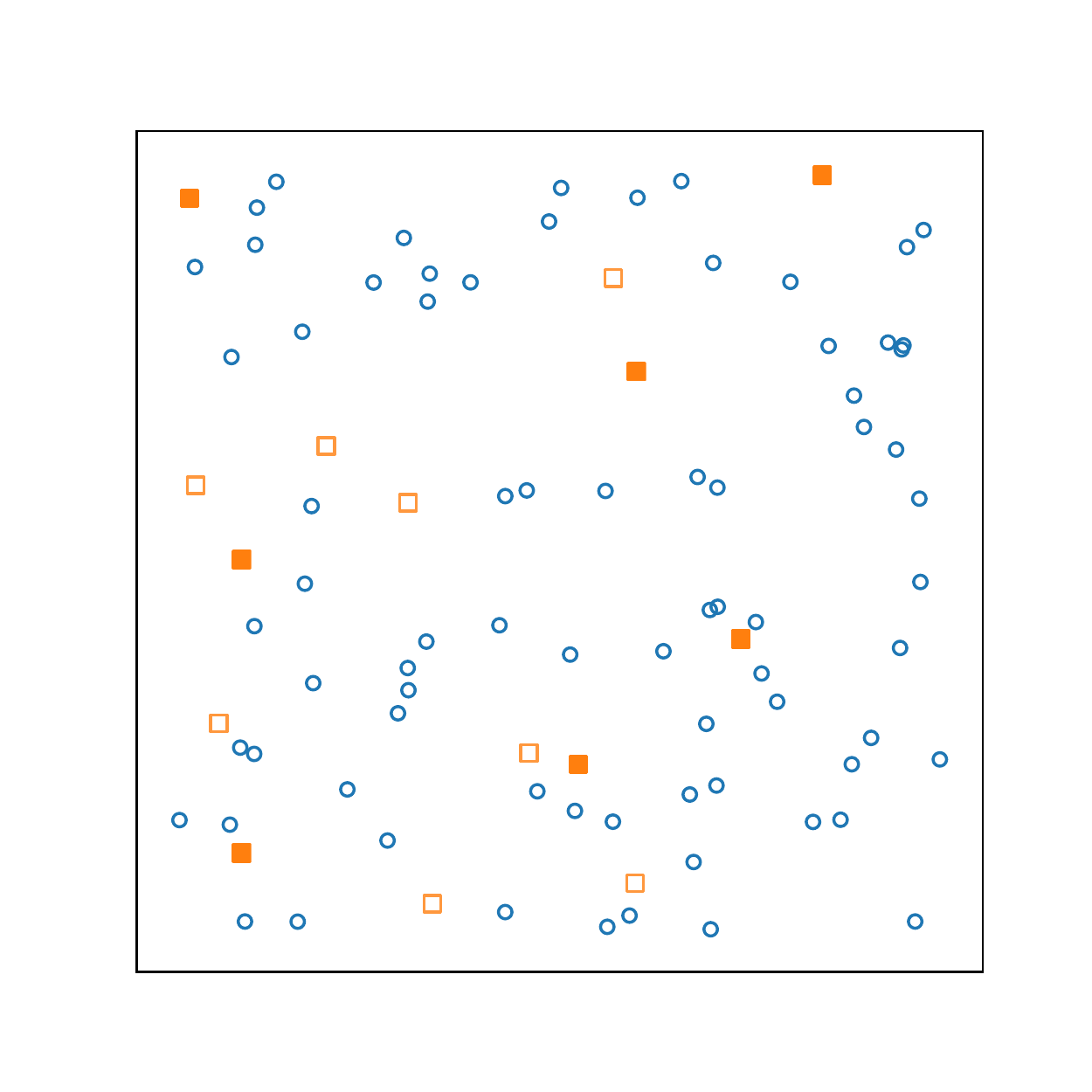}}%
	\subfloat[2nd-stage, scenario\#1]{\includegraphics[trim={1cm 1cm 1cm 1cm}, clip, scale=0.5]{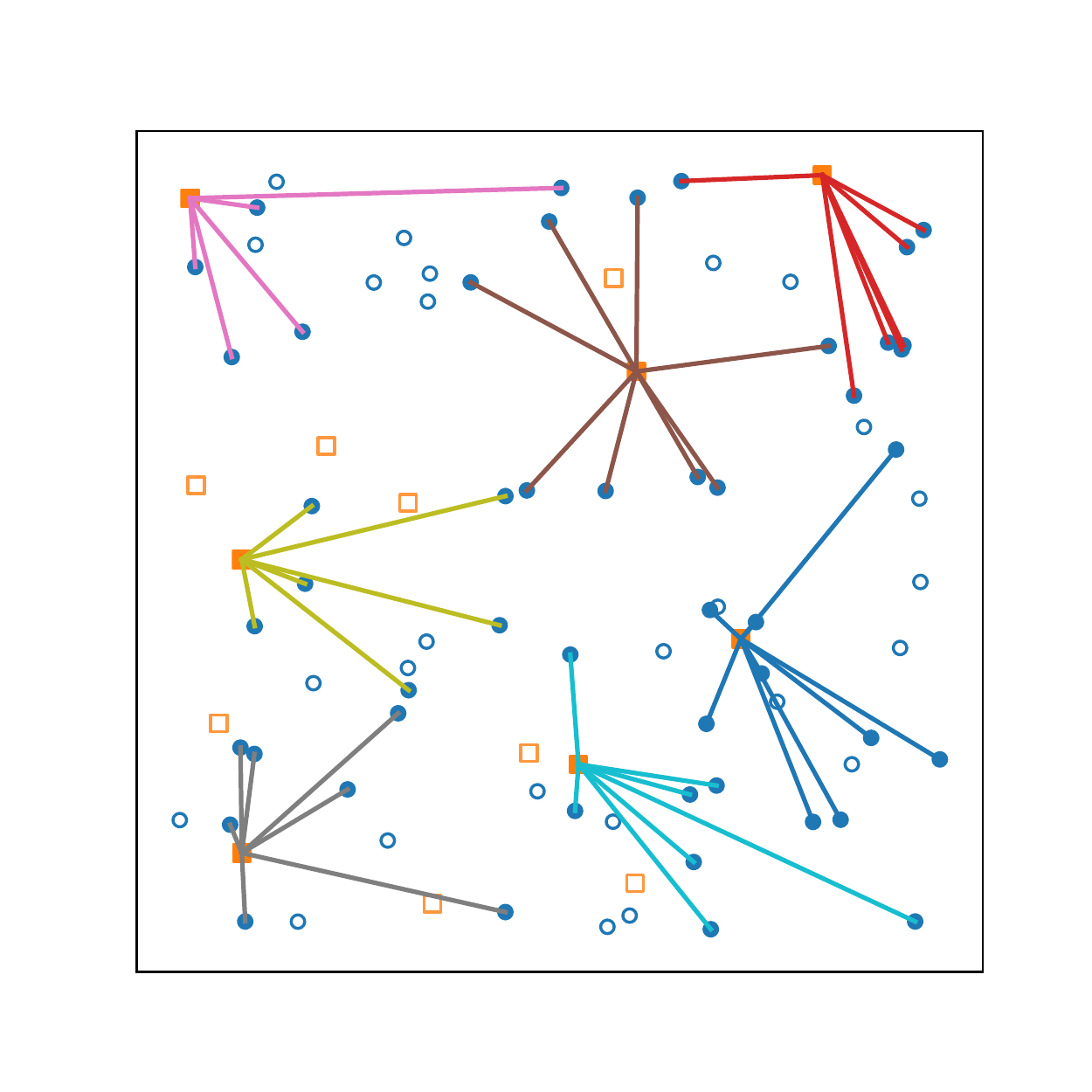}}%
	\subfloat[2nd-stage, scenario\#2]{\includegraphics[trim={1cm 1cm 1cm 1cm}, clip, scale=0.5]{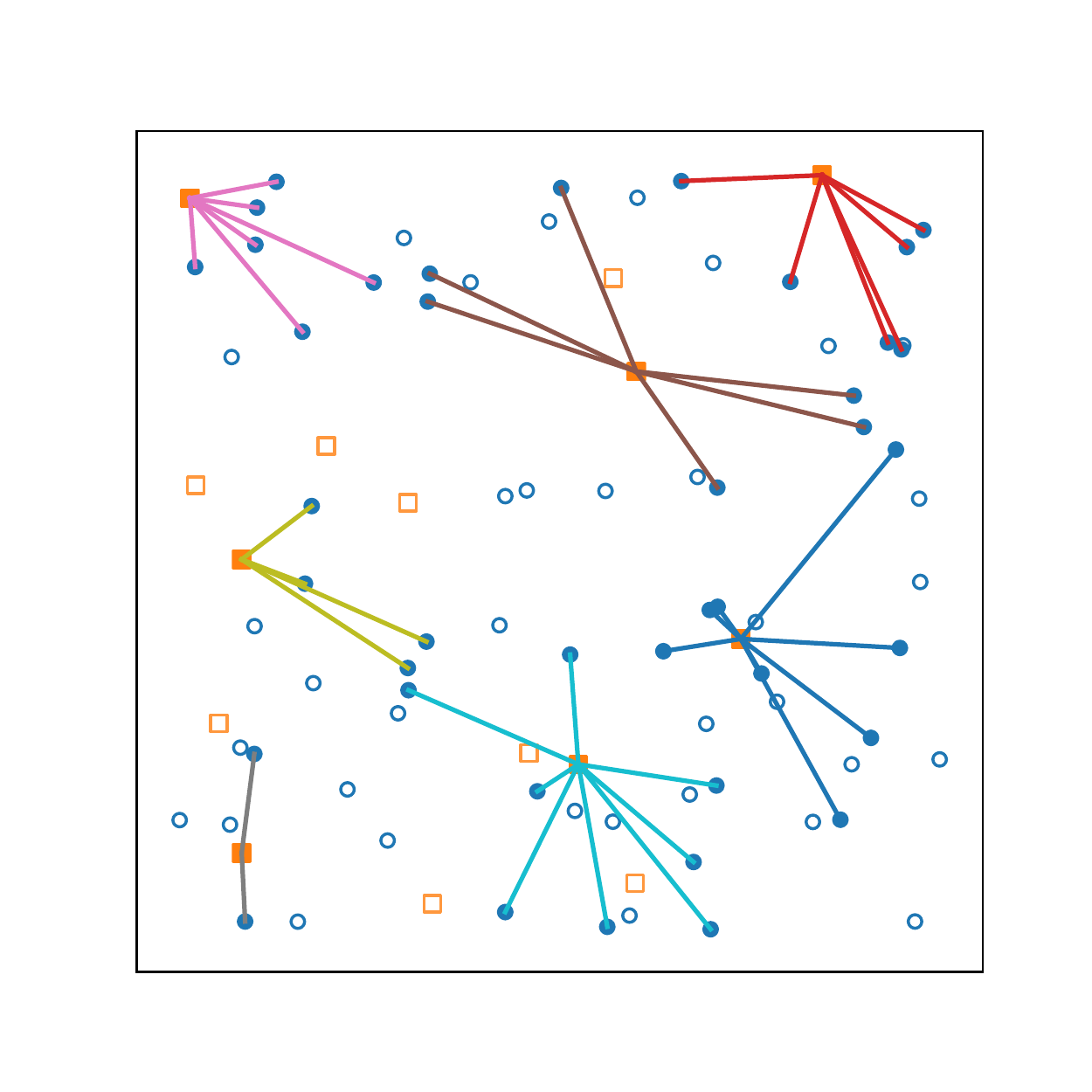}}%
	\caption{Illustration of the first-stage decisions (location of facilities) and second-stage decisions (assignment to randomly present customers) for a facility location problem with 15 facilities (squares) and 75 customers (circles).}%
	\label{FLP}%
\end{figure}

The mixed integer problem has binary first and second stage decisions. 
We denote by $F$ and $C$ the sets of predefined positions for the facilities and customers, respectively, which are represented by squares and circles in Figure~\ref{FLP}. 
The first stage variable $x_f$ can take two values: 1 if a facility is chosen to be located at position $f\in F$ and 0 otherwise. 
Located facilities are represented by a filled squares on the figure, while hollow squares are positions left empty. 
The location of a facility induces a cost $c_f$. Each facility have up to $u$ units of resources that are used to provide service to the customers.
Also, we consider that at most $v$ facilities can be opened overall.

The presence of a customer at position $c$ in scenario $i$ is given by a binary random parameter $h_{c}^{i}$.
When present, customer $c$ consumes $d_{cf}$ units of resources from the facility $f$ to which it is connected.
Each occupied customer position is represented by a filled circle in panels (b)-(c) in Figure~\ref{FLP}. 
The second stage binary variable, $y_{cf}^{i}$, is the scenario-dependent decision of matching the facility at position $f$ with the customer at position $c$ in scenario $i$. 
The cost incurred by this matching is $q_{cf}$. 
These matching decisions are represented by the links between squares and circles in panels (b)-(c) in Figure~\ref{FLP}.
Finally, if the total resources of facility $f$ is insufficient to provide service to all the customers $c$ matched to it, then $z_{f}^{i}$ defines the scenario-dependent decision representing the additional amount of capacity that is required.
In this case, a unit cost of $b_f$ is imposed for each additional unit of capacity.

The two-stage problem is formulated as follows:
\begin{align}
\label{FLP-obj}
\min & \ \sum_{f\in F} c_f x_f + \frac{1}{N} \sum_{i=1}^{N} \left(\sum_{c\in C} \sum_{f \in F} q_{cf} y_{cf}^{i} + \sum_{f\in F} b_f z_f^{i}\right) \\
\text{s.t.} & \ \sum_{f\in F} x_f \leq v \label{FLP-c1} \\
& \ \sum_{c\in C} d_{cf} y_{cf}^{i} \leq u x_f + z_f^{i} & \forall (f,i) \in F\times\setz1N \label{FLP-c2}\\
& \  z_f^{i} \leq M x_f & \forall (f,i) \in F\times\setz1N \label{FLP-c3}\\
& \  \sum_{f \in F} y_{cf}^{i} = h_{c}^{i} & \forall (c,i) \in C\times\setz1N \label{FLP-c4}\\
& \ x_f\in\{0,1\}, \ y_{cf}^{i} \in \{0,1\}, \ z_{f}^{i} \in [0,\infty) & \forall (c,f,i) \in C\times F\times\setz1N. \label{FLP-v}
\end{align}

The objective function (\ref{FLP-obj}) is composed of three terms, which from left to right are: the cost of locating facilities, the cost of matching facilities to customers, and the penalty for exceeding the facilities capacity. 
The last two are averaged over the $N$ scenarios since the decisions are scenario-dependent (but not the parameters). 
As for the constraints, (\ref{FLP-c1}) enforces the upper bound on the total number of facilities being located. 
Constraints (\ref{FLP-c2}) measure the deviations from the capacities of the facilities $f \in F$ in the scenarios $i \in \setz1N$ through the slack variables $z_{f}^{i}$, which are then penalized in the objective function. 
Constraints (\ref{FLP-c3}) guarantee that, for each facility location $f \in F$, the slack variables are zero if the facility is not located (via a big-$M$ constraint), which are necessary to ensure that no customer will be matched to it. 
To enforce that exactly one facility is assigned to a customer if the latter is present in scenario $i \in \setz1N$, constraints (\ref{FLP-c4}) are included in the model.
Finally, constraints (\ref{FLP-v}) impose the necessary integrality and non-negativity requirements on the decision variables.


We consider 10 randomly-generated instances of the problem where each instance has 100 and 1000 possible locations for the facilities and the clients, respectively. These locations are independently sampled from a uniform distribution in $[0,1]^2$. Each instance has a set of $100$ scenarios of client presence, which are sampled in two steps: first, $100$ random numbers $p\in(0.2,0.8)$ are independently sampled from a uniform distribution, then each number is used as the parameter of a Bernoulli distribution $\mathcal{B}(p)$ from which the presence of the 1000 clients is sampled. In other words, we have $h_{c}^{i} \sim \mathcal{B}(p_i)$ for each $c\in C$ with $p_i\sim U(0.2,0.8)$ for all $i\in\setz1N$. The purpose of randomly generating the Bernoulli parameter $p$, as opposed to setting a fixed one, is to increase the variability of customers across the scenarios, since for 1000 possible locations a fixed parameter $p$ will typically results in all scenarios having roughly the same number of customers (about 1000$p$), while a uniformly sampled $p\in(0.2,0.8)$ will produce scenarios with a range of about 200 to 800 customers.

For this number of scenarios (100), customers (1000), and facilities (100), each instance of the problem has about 10 million binary variables. This makes the solution of the original mixed-integer problem highly challenging even for a top-notch commercial solver, which on the full problem with 100 scenarios could only reach a gap of 8\% to 15\% within 6 hours of computation (see left column of Table~\ref{Table_CPLEX_vs_CSSC_gap}). For most of the instances (7 out of 10), the computation even stopped before the time limit due to a lack of memory in the branch-and-bound process (see left column of Table~\ref{Table_CPLEX_vs_CSSC_time}). 

To overcome this limitation, we consider solving the problem in a deterministic form, i.e., using a single scenario provided by the following scenario reduction methods: CSSC, $k$-medoids, $k$-modes, and Monte Carlo. All these methods have in common that the scenario they output remains in a binary format, which is necessary to have a sensible facility location problem (a fractional presence of customers in constraint \eqref{FLP-c3} will lead to an unfeasible problem). The deterministic one-scenario problem is then solved twice with different time limits (300 and 600 seconds) and the solution is evaluated in the original problem with 100 scenarios to assess its true value. The solution obtained in the best of the two runs is the one displayed in Table~\ref{Table_CPLEX_vs_CSSC_gap}, where the gaps (which by definition range from 0\% to 100\%) are computed with respect to the bound that was obtained when solving the original problem. 

We see that the CSSC algorithm is the only method that provides solutions with consistently low error (less than 3\%), which is significantly less than the error made after 6 hours of solving the original problem. By comparison, $k$-medoids and Monte Carlo manage to find a good quality solution for only one and two instances out of 10, respectively, and $k$-modes could not output a single solution with a gap less than 90\%. This extremely bad results of $k$-modes could be explained by the fact that it is the only method outputting scenarios that are not in the original set.
 
The computational times to run the CSSC algorithm are very modest as compared to the time spent solving the original problem, as seen in Table~\ref{Table_CPLEX_vs_CSSC_time} where we display the longest of the two runs. In fact, it takes generally less than 30 min for the algorithm to go through all the different steps: computation of the opportunity cost matrix (column ``step 1''), solution of the MIP partitioning problem (column ``step 2''), solution of the approximate problem (column ``approx.''), and evaluation of the solution in the original problem (column ``eval.''). 
We also stress that the computation of the opportunity cost matrix was done on the linear relaxation of the problem, but without parallelization, so the times could be further improved.

\begin{figure}[h]
	\centering
	\includegraphics[trim={0 0 0 0}, clip, scale=0.7]{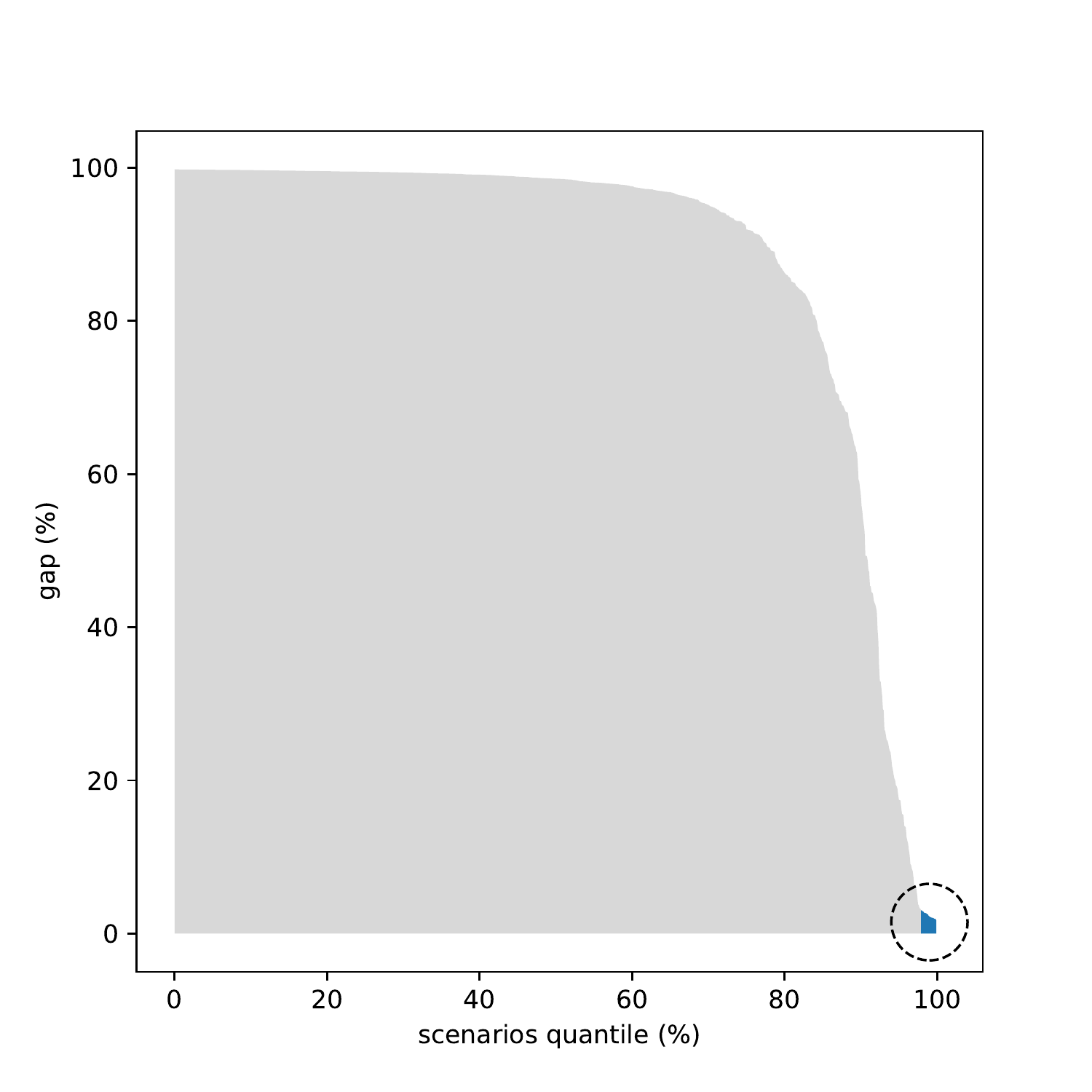}
	\caption{Distribution of optimality gaps across the 100 scenarios and 10 instances of the facility location problem. All scenarios provided by the CSSC algorithm lie in the bottom-right blue part of the distribution.}%
	\label{Dist-gap}%
\end{figure}

Finally, to evaluate more precisely the ability of the CSSC method to single out what appears to be some of the best scenarios out of the batch of 100 of them, we have solved the 10 instances on each one of the 100 scenarios individually (hence a total of 1000 one-scenario problems). This way we want to make sure that we are not facing a situation where the majority of scenarios would be appropriate to build the approximate problem. In fact, the results displayed in Figure~\ref{Dist-gap} show the opposite: the vast majority of scenarios are completely unsuitable to approximate the original problem, as only a few percent of the scenarios (about 2-3\%) lead to a problem with less than 5\% error gap (see the circle at the bottom right of the figure). The CSSC algorithm was able to always select one of them. More specifically, over the ten instances it selected the best scenario six times, the second best twice and the third best twice as well. 

\begin{table}[h]
	\centering
	\begin{tabular}{|c|c|cccc|}
		\hline		
		\multicolumn{1}{|c|}{} & \multicolumn{1}{c|}{} & \multicolumn{4}{c|}{$K=1$} \\
		\hline
		instance & original problem & CSSC & $k$-medoids & $k$-modes & MC \\
		\hline
		0 & 13.2 & \textbf{2.2} & 99.7 & 99.3 & 97.7 \\
		1 & 7.7 & \textbf{2.9} & 99.7 & 99.5 & 98.8 \\
		2 & 17.2 & \textbf{2.6} & 99.7 & 99.4 & 98.8 \\
		3 & 14.3 & \textbf{1.9} & 99.7 & 99.4 & 99.4 \\
		4 & 14.2 & 2.4 & \textbf{1.9} & 91.2 & 71.3 \\
		5 & 11.0 & \textbf{1.8} & 21.8 & 97.8 & 2.4 \\
		6 & 14.3 & \textbf{2.4} & 99.7 & 99.6 & 24.0 \\
		7 & 11.6 & \textbf{1.8} & 44.8 & 97.2 & 97.1 \\
		8 & 11.2 & \textbf{1.7} & 79.2 & 98.4 & 2.0 \\
		9 & 9.6 & \textbf{2.1} & 95.8 & 98.9 & 84.9 \\
		\hline
	\end{tabular}
	\caption{Comparison between the gap (in percentage) reached by solving the original problem with 100 scenarios (with a time limit of 6 hours) and the gaps of the solutions with a single scenario ($K=1$) given by the CSSC algorithm, $k$-medoids, $k$-modes, and Monte Carlo. (Data in bold font singles out the smallest value row-wise.)}
	\label{Table_CPLEX_vs_CSSC_gap}
\end{table}

\begin{table}[h]
	\centering
	\begin{tabular}{|c|c|cccc|c|}
		\hline
		\multicolumn{1}{|c|}{} & \multicolumn{1}{c|}{} & \multicolumn{5}{c|}{CSSC ($K=1$)} \\
		\hline
		instance & original problem & step 1 & step 2 & approx. & eval. & total  \\
		\hline
		0 & 17405* & 1098 & 5.2 & 601 & 49 & 1754 \\ 
		1 & 21610 & 1207 & 5.2 & 602 & 82 & 1896 \\ 
		2 & 21605 & 1216 & 5.1 & 601 & 113 & 1935 \\ 
		3 & 15028* & 1213 & 5.2 & 601 & 115 & 1935 \\ 
		4 & 20047* & 1317 & 5.2 & 602 & 115 & 2039 \\ 
		5 & 21607 & 1322 & 5.2 & 601 & 399 & 2327 \\ 
		6 & 13172* & 1263 & 5.3 & 602 & 200 & 2071 \\ 
		7 & 17303* & 1334 & 5.3 & 602 & 242 & 2183 \\ 
		8 & 20297* & 1368 & 5.2 & 601 & 193 & 2166 \\ 
		9 & 20033* & 1379 & 5.1 & 601 & 199 & 2184 \\ 
		\hline
	\end{tabular}
	\caption{Comparison between the computational times (in seconds) of the solution of the original problem with 100 scenarios (limited to 6 hours) and those of the different steps of the CSSC algorithm. (An asterisk (*) signifies that the solution of the original problem stopped before the time limit with an out-of-memory status.)}
	\label{Table_CPLEX_vs_CSSC_time}
\end{table}

\section{Conclusion}
\label{sec:conclusion}

We have developed a new scenario reduction algorithm that can select appropriately a subset of scenarios when the original stochastic problem is too large to be solved using all of them. This algorithm works directly with the scenario data (no underlying probability distribution is needed), and it uses the problem itself to define a measure of proximity between scenarios. This means that no arbitrary metric in the space of random outcomes must be specified; it is actually the cost function of the problem that provides this metric, which is then used to cluster the scenarios together, select some representatives, and assign them probabilities. For this reason the clustering algorithm is said to work in the space of cost values ($\R$) as it has no view of the original space ($\R^d$) in which the scenarios are expressed. 

We have tested this new algorithm on two stochastic mixed-integer problems: a network design and a facility location problem. In both cases, the cost space scenario clustering (CSSC) algorithm largely outperforms the other methods that cluster scenarios based on their proximity in the space of random outcomes. For instance, in the network design problem, results show that with only 10 scenarios out of 75, the CSSC algorithm solves more than 85\% of instances with less than 10\% error, as opposed to a mere 36\% for the second best tested method, and only 1\% for random selection via Monte Carlo sampling. For the facility location problem, all instances are better solved with a single scenario provided by the CSSC algorithm (in less than 40 min) than for the actual original problem with 100 scenarios solved for 6 hours. 

These better performance results should be traded-off against a larger computational cost induced by the CSSC algorithm as compared to other clustering methods that do not account for the problem. One direction of future research goes towards reducing this additional cost. For instance, the computation of the opportunity cost matrix could be done using a relaxation of the problem (linear or other relaxations). Also, the solution of one problem could be used as an initial start for the other problems to be solved, as problems may only slightly change from one iteration to the next. As for the mixed-integer partitioning problem, it could be solved using heuristic methods rather than to global optimality as we did in this paper. Finally, another direction of research would be to generalize the clustering criterion to multistage problems.

\FloatBarrier
\newpage
\appendix

\bibliographystyle{apalike}
\bibliography{References}

\clearpage

\end{document}